\documentclass[a4paper, reqno]{amsart}

\usepackage{amsmath,amsfonts,amssymb,amsthm,graphicx,mathtools}
\usepackage{mathrsfs}   % Script fonts
\usepackage{xspace}     % For handling spaces after our defined commands
\usepackage{hyperref}   % Optional: for hyperlinked references
\DeclareMathOperator{\aut}{Aut}
\DeclareMathOperator{\Aut}{Aut}
\DeclareMathOperator{\AAut}{AAut}
\DeclareMathOperator{\Sym}{Sym}
\DeclareMathOperator{\sym}{Sym}
\newcommand{\N}{\mathbb{N}}
\providecommand{\Wr}{\mathop{\rm Wr}\nolimits}
\providecommand{\lcol}{\mathcal{L}}
\newcommand{\Z}{\mathbb{Z}}
\newcommand{\U}{\mathbf{U}}
\newcommand{\inv}{^{-1}}
\newcommand{\G}{\mathcal{G}}
\newcommand{\propP}[1]{(\mathrm{P}_{#1})}
\newcommand{\tdlc}{t.d.l.c.\@\xspace}
\newcommand{\scpo}{scopo\@\xspace}
\newcommand{\scpos}{scopos\@\xspace}

\usepackage[T1]{fontenc}
\usepackage{tikz}
\usetikzlibrary{positioning}
\usetikzlibrary{intersections, decorations.pathreplacing, decorations.markings}
\usetikzlibrary{calc}
\usetikzlibrary{shapes.misc}
\usepackage{xifthen}
\tikzstyle{Vertex}=[shape=circle, draw=black, fill=black, inner sep=0pt, minimum size=3]
\tikzstyle{OpenVertex}=[shape=circle, draw=black, fill=white, inner sep=0pt, minimum size=3]
\tikzstyle{RedVertex}=[shape=circle, draw=red, fill=red, inner sep=0pt, minimum size=3]
\tikzstyle{TextNode}=[shape=rectangle, inner sep=0pt, minimum size=3]
\usetikzlibrary{decorations.pathmorphing}
\tikzset{snake it/.style={decorate, decoration=snake}}
\tikzset{->-/.style={decoration={
  markings,
  mark=at position .5 with {\arrow{>}}},postaction={decorate}}}

\newtheorem{theorem}{Theorem}
\newtheorem*{theorem*}{Theorem}

\newtheorem*{lemma*}{Lemma}
\newtheorem{corollary}[theorem]{Corollary}
\newtheorem*{corollary*}{Corollary}

\newtheorem{proposition}[theorem]{Proposition}
\newtheorem*{proposition*}{Proposition}
\newtheorem{project}{Project}

\theoremstyle{definition}
\newtheorem{definition}[theorem]{Definition}
\newtheorem*{definition*}{Definition}

\newtheorem{remark}[theorem]{Remark}
\newtheorem{example}[theorem]{Example}
\newtheorem*{example*}{Example}
\newtheorem*{remark*}{Remark}
\newtheorem*{exercise*}{Exercise}
\newtheorem*{furtherprojects}{Further projects}

\begin{document}

\author{Colin D.~Reid}
\address{Colin D.~Reid. The University of Newcastle, School of Mathematical and Physical Sciences, Callaghan, NSW 2308, Australia. Email: colin@reidit.net}

\author{Simon M.~Smith}
\address{Simon M.~Smith. Charlotte Scott Research Centre for Algebra, University of Lincoln, Lincoln, U.K. Email: sismith@lincoln.ac.uk}

\title[Local action diagrams]{An introduction to the local-to-global behaviour of groups acting on trees and the theory of local action diagrams}

\begin{abstract}
The primary tool for analysing groups acting on trees is Bass--Serre Theory.
It is comprised of two parts: a decomposition result, in which an action is decomposed via a graph of groups, and a construction result, in which graphs of groups are used to build examples of groups acting on trees. The usefulness of the latter for constructing new examples of `large' (e.g.~nondiscrete) groups acting on trees is severely limited. There is a pressing need for new examples of such groups as they play an important role in the theory of locally compact groups. An alternative `local-to-global' approach to the study of groups acting on trees has recently emerged, inspired by a paper of Marc Burger and Shahar Mozes, based on groups that are `universal' with respect to some specified `local' action.
In recent work, the authors of this survey article have developed a general theory of universal groups of local actions, that
behaves, in many respects, like Bass--Serre Theory. We call this the theory of local action diagrams. The theory is powerful enough to completely describe all closed groups of automorphisms of trees that enjoy Tits' Independence Property $\propP{}$.

This article is an introductory survey of the local-to-global behaviour of groups acting on trees and the theory of local action diagrams. The article contains many ideas for future research projects.
\end{abstract}

\maketitle

\section{Introduction}
Actions on trees have a significant role in the general theory of finite and infinite groups. 
In finite group theory such actions are, for example, a natural setting for questions about vertex-transitive groups $G$ of automorphisms of connected finite graphs. For such a pair $G \leq \aut(\Gamma)$ with $G$ nontrivial, the universal cover of $\Gamma$ is an infinite regular tree $T$, and there is a natural projection $\pi$ of $T$ to $\Gamma$. 
The
  fundamental\footnote{We use the notation $\Pi(\Gamma, v)$ for the fundamental group, rather than the more common $\pi(\Gamma, v)$, because for us $\pi$ will be reserved for projection maps.  } group  
$\Pi(\Gamma, v)$ of $\Gamma$ at any vertex $v$ can be identified with a subgroup of $\aut(T)$ in such a way that $\Gamma$ can be identified with the quotient graph $\Pi(\Gamma, v) \backslash T$ and the lift $\tilde{G}$ of $G$ along $\pi$ can be identified with the quotient $\tilde{G} / \Pi(\Gamma, v)$. For a thorough description of this correspondence, see \cite{potocnik_spiga_19} for example.

In this context, important open problems in finite groups can translate into natural questions about groups acting on trees. Consider, for example, the well-known Weiss conjecture (\cite{weiss78}) due to Richard Weiss. In the language of finite groups, it states that there exists a function $f: \N \rightarrow \N$ such that if $G$ is vertex-transitive and locally primitive\footnote{For a group $G$ acting on a graph $\Gamma$, the {\it neighbours} of a vertex $v$ are the vertices in $\Gamma$ at distance one from $v$. Each vertex stabiliser $G_v$ induces a permutation group on the neighbours of $v$; if every such induced permutation group has some permutational property $\mathcal{P}$ (e.g. transitive, primitive), we say that $G$ is {\it locally $\mathcal{P}$}.} on a graph with finite valency $k$, then all vertex stabilisers satisfy $|G_v| < f(k)$. Equivalently, in the language of groups acting on trees, it states that for each $k \in \N_{\geq 3}$ the automorphism group of the $k$-regular tree $T_k$ contains only finitely many conjugacy classes of discrete, locally primitive and vertex transitive subgroups.

The role played by actions on trees in the general theory of infinite groups is, of course, well-known. Much of our understanding of groups acting on trees comes from the celebrated 
Bass--Serre Theory, described in detail in Jean-Pierre Serre's book \cite{Serre:trees}. We give an introductory overview to Bass--Serre Theory in Section~\ref{Sec:BassSerre}.
  Serre's theory concerns a group $G$ {\it acting on a tree} $T$; by this we mean that $G$ acts on $T$ as a group  of automorphisms of $T$. We will denote such an action by the pair $(T,G)$, and will often identify $(T,G)$ with its image in the group $\aut(T)$ of automorphisms of $T$. In Bass--Serre theory,  
the algebraic structure of a group $G$ acting on a tree $T$ is `decomposed' into pieces with the decomposition described via a combinatorial structure called a {\em graph of groups}. This graph of groups associated to the action $(T,G)$ is essentially a vertex- and edge-coloured graph, where each colour is in fact a group. The groups in this palette are called the {\it vertex groups} and {\it edge groups} of the graph of groups. Importantly, this decomposition process can be reversed, where one starts with a graph of groups $\Gamma$ and then constructs the universal cover $\tilde{T}$ (which is a tree) of $\Gamma$ and the fundamental group $\Pi$ of $\Gamma$. The group $\Pi$ acts on $\tilde{T}$ in such a way that its associated graph of groups is again $\Gamma$. This `decomposition' of $(T,G)$ in Bass--Serre Theory hinges on the observation that if $\Gamma$ is the associated graph of groups for $(T,G)$, and $\tilde{T}$ and $\Pi$ are respectively the universal cover and fundamental group of $\Gamma$, then the actions $(T,G)$ and $(\tilde{T}, \Pi)$ can be identified.

Two algebraic constructions, the HNN extension and the amalgamated free product, play a foundational role in this theory: the graph of groups of an HNN extension is a single vertex with a loop, and the graph of groups of an amalgamated free product is a pair of vertices joined by an edge. Intuitively, edges and loops in the graph of groups of the action give information about how the action decomposes into HNN extensions and amalgamated free products of the vertex groups. 

Naturally, there are areas in which the usefulness of Bass--Serre Theory is limited. One of the most significant limitations can be seen when attempting
 to construct a group acting on a tree with certain desired properties, by first writing down a graph of groups. We discuss this in Section~\ref{subsection:bass_serre_limitations}.   Essentially   Bass--Serre Theory can be used to construct an action on a tree with a given graph of groups, but the action itself cannot, in general, be fully controlled. This limitation is the source of the intractability of several famous conjectures, including the aforementioned Weiss Conjecture and the related Goldschmidt--Sims Conjecture (see \cite{BurgerMozes}, \cite{goldschmidt}). This limitation is particularly problematic when trying to construct nondiscrete actions on trees, and recent developments in the theory of locally compact groups have made the need to address this limitation more acute.

In the theory of locally compact groups, the structure theory of compactly generated locally compact groups is known to depend on the class $\mathscr{S}$ of nondiscrete, compactly generated, locally compact topologically simple groups; see Pierre-Emmanuel Caprace and Nicolas Monod's paper \cite{CapraceMonod}. For a survey of the class $\mathscr{S}$, placing it in a broad mathematical context, see Caprace's survey paper \cite{CapraceSimple}.
Groups acting on trees are one of the main sources of examples of groups in $\mathscr{S}$. In fact the very first examples of nonlinear, nondiscrete, locally compact simple groups were constructed using groups acting on trees, by Jacques Tits (answering a question of J.~P.~Serre) in \cite{Tits70}. Tits used an independence condition he called Property $\propP{}$ for groups acting on trees, and showed that for any group $G$ of automorphisms of a tree $T$, if $G$ has Property $\propP{}$, fixes (setwise) no nonempty proper subtree of $T$ and fixes no end of $T$, then the subgroup $G^+$ of $G$ generated by arc stabilisers is (abstractly) simple\footnote{
% Footnote
We will always think of $\Aut(T)$ as a topological group under the {\it permutation topology}, defined below.
A topological group is {\it topologically simple} if it has no nontrivial proper closed normal subgroups; it is {\it abstractly simple} if it has no nontrivial proper normal subgroups. When we write simple, we will always mean abstractly simple.
}
(see Section~\ref{Section:tits}). As an immediate consequence, for the infinite regular tree $T_n$ of (finite or infinite) valency $n$ we have that $\aut(T_n)^+$ is simple for $n \geq 3$. This group is nondiscrete and when $n$ is finite it is easily seen to be compactly generated and locally compact (under the permutation topology).

The limitations of Bass--Serre Theory when constructing nondiscrete groups acting on trees can be avoided using a variety of techniques inspired by a 2000 article \cite{BurgerMozes} by Marc Burger and Shahar Mozes; these techniques might collectively be known as {\it local-to-global constructions of groups acting on trees}. In their paper, Burger and Mozes take a permutation group $F \leq S_n$ of finite degree $n$ and build a group $\U(F)$ of automorphisms of $T_n$ that has {\it local action} $F$ (or is {\it locally-$F$}); that is, vertex stabilisers in $\U(F)$ induce $F$ on the set of neighbours of the stabilised vertex. For example, $\U(S_n) = \aut(T_n)$. These {\it Burger--Mozes} groups $\U(F)$ enjoy Tits' Independence Property $\propP{}$, and when $F$ is transitive $\U(F)$ is `universal' among subgroups of $\Aut(T_n)$ that have local action $F$; that is, $\U(F)$ contains an $\Aut(T_n)$-conjugate of every subgroup of $\Aut(T_n)$ that is locally-$F$. When $F$ is transitive and generated by its point stabilisers, $\U(F)$ contains a simple subgroup $\U(F)^+$ of index $2$. 
Since 2000, the majority of constructions of compactly generated, simple, locally compact groups have used the ideas of Tits \cite{Tits70} and Burger--Mozes \cite{BurgerMozes}.

In \cite{SmithDuke}, the second author generalised the Burger--Mozes construction to {\it biregular trees}   $T_{m,n}$,   where $m,n \geq 2$ are (possibly infinite) cardinals and $T_{m,n}$ is the infinite tree in which vertices in one part of its bipartition have valency $m$ and those in the other part have valency $n$.
Given permutation groups $F_1 \leq S_m$ and $F_2 \leq S_n$,
this generalisation is a group $\U(F_1, F_2) \leq \aut(T_{m,n})$, called {\it the box product} of $F_1$ and $F_2$. The group $\U(F_1, F_2)$ has local actions $F_1$ (at vertices of $T_{m,n}$ in one part of the bipartition) and $F_2$ (at vertices in the other part of the bipartition); we say any such action on $T_{m,n}$ is {\it locally-$(F_1, F_2)$}. 
For $m \geq 3$ and $F \leq S_m$ it can be shown that the Burger--Mozes group $\U(F)$ is isomorphic as a topological group (but not as a permutation group) to $\U(F, S_2)$. The properties of $\U(F_1,  F_2)$ mirror those of the Burger--Mozes groups: it has Tits' Independence Property $\propP{}$, and when $F_1$ and $F_2$ are transitive, $\U(F_1, F_2)$ is `universal' among subgroups of $\Aut(T_{m,n})$ that are locally-$(F_1, F_2)$, in that it contains an $\aut(T_{m,n})$-conjugate of every locally-$(F_1,F_2)$ subgroup of $\aut(T_{m,n})$. When $F_1$ and $F_2$ are generated by their point stabilisers and at least one group is nontrivial, $\U(F_1,F_2)$ is simple if and only if $F_1$ or $F_2$   is   transitive.

The Burger--Mozes groups are automatically totally disconnected, locally compact (henceforth, ``\tdlc'') and compactly generated groups because they   are   defined only on locally-finite trees. Despite admitting actions on non-locally finite trees, $\U(F_1, F_2)$ can still be a compactly generated \tdlc group under mild conditions on $F_1$ and $F_2$. Indeed, if $F_1$ and $F_2$ are closed and compactly generated, with compact nontrivial point stabilisers and finitely many orbits, and either $F_1$ or $F_2$ is transitive (e.g. $F_1$ is closed, nonregular, subdegree-finite and primitive and $F_2$ is finite, nonregular and primitive), then $\U(F_1, F_2)$ is a nondiscrete compactly generated \tdlc group. By taking $F_1$ to be infinite permutation representations of various   Tarski--Ol'shanski{\u\i}   Monsters (see \cite{olshanski}), the second author used this box product construction to show that there are precisely $2^{\aleph_0}$ isomorphism types of nondiscrete, compactly generated simple locally compact groups, answering a well-known open question; the analogous result for discrete groups was proved in 1953 by Ruth Camm in \cite{Cam53}.\\

In light of the Burger--Mozes and box product constructions, the authors of this article wished to create a unified way of seeing $\U(F)$ and $\U(F_1, F_2)$, within a framework permitting yet more local actions to be specified. What emerged from this endeavour was something far deeper: a general theory of `universal' groups of local actions, that behaves, in many respects, like Bass--Serre Theory. This work, which we call {\it the theory of local action diagrams}, is fully described in our paper \cite{ReidSmith}.

To better understand this theory, let us define the {\em $\propP{}$-closure} of an action $(T, G)$ of a group $G$ on a tree $T$ as being the smallest closed subgroup of $\aut(T)$ with Tits Independence Property $\propP{}$ that contains
  $(T,G)$; we denote it $G^{\propP{}}$. That is, $G^{\propP{}}$ is the smallest closed subgroup of $\aut(T)$ with Tits Independence Property $\propP{}$ that contains the subgroup of $\aut(T)$ induced by the action of $G$ on $T$. If $G = G^{\propP{}}$ we say that $G$ is {\it $\propP{}$-closed}. \label{def:P-closure}  
In the theory, any group $G$ acting on a tree $T$ is `decomposed' into its local actions with the decomposition described via a combinatorial structure called a {\it local action diagram}. This local action diagram is a graph decorated with sets and groups that codify the local actions of $G$. 
Importantly, this decomposition process can be reversed, where one starts with a local action diagram $\Delta$ and then constructs an arc-coloured tree ${\bf T}$ called the $\Delta$-tree and a group $\U(\Delta)$ called the {\it universal group of $\Delta$}. The group $\U(\Delta)$ acts on ${\bf T}$ in such a way that its local action diagram is again $\Delta$, and moreover it exhibits various desirable global properties that are often impossible to verify for groups arising from graphs of groups via Bass--Serre Theory. This `decomposition' of $(T,G)$ into local actions hinges on the observation that if $\Delta$ is the associated local action diagram for $(T,G)$ and ${\bf T}$ and $\U(\Delta)$ are the $\Delta$-tree and the universal group of $\Delta$, then the actions $(T,G^{\propP{}})$ and $({\bf T}, \U(\Delta))$ can be identified. Thus, viewing $\U(\Delta)$ as a subgroup of $\aut(T)$, we see that $\U(\Delta)$ is `universal' with respect to the local actions of $G$; that is, $\U(\Delta)$ contains an $\aut(T)$-conjugate of every action $(T,H)$ whose local action diagram is also $\Delta$.

The Burger--Mozes groups $\U(F)$ and the box product construction $\U(F_1, F_2)$ play a foundational role in this new theory: the local action diagram of $\U(F)$ is a (suitably decorated) graph consisting of a single vertex with a set of loops, each of which is its own reverse; the local action diagram of $\U(F_1, F_2)$ when $F_1$ and $F_2$ are transitive is a (suitably decorated) graph consisting of a pair of vertices joined by an edge and no loops (c.f.~the graphs of groups of HNN extensions and amalgamated free products).

The theory of local action diagrams is powerful enough to give a complete description of closed actions on trees with Tits' Independence Property $\propP{}$: these actions are precisely the universal groups of local action diagrams. From this one immediately obtains a robust classification of all actions $(T,G)$ with Tits' Independence Property $\propP{}$, by taking closures in the permutation topology. Note that $(T,G)$ and its closure in the permutation topology have the same orbits on all finite (ordered and unordered) sets of vertices.

Working with local action diagrams is remarkably easy; unlike for graphs of groups there are no embeddability issues to contend with. Moreover, all properties of faithful closed actions with property $\propP{}$ can be read directly from the local action diagram; many can be read easily from the diagram, including --- surprisingly ---  geometric density, compact generation and simplicity.\\

The theory of local action diagrams is described in \cite{ReidSmith}. This note is intended to be an accessible introduction to the theory, placing it in its broad context and giving illustrative diagrams and examples, while omitting most proofs. It largely follows the contents of the second   author's   plenary talk at Groups St Andrews 2022 in Newcastle.
Our intention with this note is to make the motivations and ideas in the theory accessible to non-specialists. To that end, we give an overview of Bass--Serre Theory, with enough depth so that the reader can understand its limitations, appreciate how the theory of local action diagrams mirrors the rich relationship between action and combinatorial description in Bass--Serre Theory, and to see how fundamentally different the two theories are.
Following this we introduce the theory of local action diagrams, giving examples and some consequences, but largely omitting proofs (all proofs can be found in \cite{ReidSmith}). We conclude with some suggestions for future research projects.

Since this note is intended to be an introduction to local action diagrams, many topics from \cite{ReidSmith} have been omitted. The most significant omissions concern subgroups of universal groups of local action diagrams. We refer the interested reader to \cite{ReidSmith} for further details.

\subsection*{Notation and conventions}
Unless otherwise stated we follow the definitions in our paper \cite{ReidSmith}. In particular, our graphs can have multiple distinct edges between two vertices, and each edge is comprised of two arcs (one in each direction). Loops are allowed, meaning that our graphs are graphs in the sense of Serre (except for us a loop $a$ may or may not equal its reverse $\overline{a}$, which is not the case for Serre). Our graphs $\Gamma$ have a   (finite or infinite)   vertex set $V = V\Gamma$, a   (finite or infinite)   arc set $A = A\Gamma$, an arc-reversal map $a \mapsto \overline{a}$ (also called an arc or edge {\it inversion}) together with an {\it origin} map $o:A \rightarrow V$ and a {\it terminal} map $t:A \rightarrow V$, so that $a \in A$ is an arc {\it from $o(a)$ to $t(a)$}. Edges are pairs $\{a, \overline{a}\}$ and are said to {\it contain} vertices $o(a)$ and $t(a)$, or to be {\it between} $o(a)$ and $t(a)$. We define {\it loops} to be arcs $a$ such that $o(a) = t(a)$. 
The {\it valency} of a vertex $v$ is $|o\inv(v)|$, sometimes denoted $|v|$; if this is finite for all vertices then $\Gamma$ is {\it locally finite}.
A {\it leaf} in $\Gamma$ is a vertex with exactly one edge containing it, and that edge is not a loop. A graph is {\it simple} if it has no loops and there is at most one edge between any two vertices.

Our graphs are not simple so we must define paths with care. Given an interval $I \subseteq \Z$, let $\hat{I} = \{i \in I : i+1 \in I\}$; a {\it path} in $\Gamma$ indexed by $I$ is then a sequence of vertices $(v_i)_{i \in I}$ and edges $(\{a_i, \overline{a_i}\})_{i \in \hat{I}}$ such that $\{a_i, \overline{a_i}\}$ is an edge in $\Gamma$ between $v_i$ and $v_{i+1}$ for all $i \in \hat{I}$. For finite $I$ the path has {\it length} $|\hat{I}|$. A path is {\it simple} if all its vertices $v_i$ are distinct from one another. We can now define {\it directed paths} in the obvious way, as a sequence of vertices and arcs. For $n > 0$, if $I = \{0,\dots,n\}$ and $v_0 = v_{n}$ and vertices $\{v_0, \dots, v_{n-1}\}$ are distinct, then the path is called a {\it cycle} of length $n$. The {\it distance} between two vertices $v,w$, denoted $d(v,w)$, is the length of the shortest path between them if it exists, and is infinite otherwise.
A graph is {\it connected} if there is a path between any two distinct vertices. For a vertex $v$ the set of vertices whose distance from $v$ is at most $k$ is called a {\it $k$-ball} and is denoted $B_v(k)$. We will sometimes write $B(v)$ for the set $B_v(1)$ of {\it neighbours} of $v$.
An {\it orientation} of $\Gamma$ is a subset $O \subseteq A\Gamma$ such that for each $a \in A\Gamma$, either $a$ or $\overline{a}$ is in $O$, but not both. For graphs that are not simple, the graph subtraction operation is not well behaved and so we avoid it except for the following situation: for $a \in A\Gamma$ the graph $\Gamma \setminus \{a\}$ is obtained from $\Gamma$ by removing arcs $a, \overline{a}$. For a simple graph $\Gamma$ with subgraph $\Lambda$ we define graph subtraction $\Gamma \setminus \Lambda$ in the usual 
  way: $\Gamma \setminus \Lambda$ is obtained from $\Gamma$ by removing all vertices that lie in $\Lambda$ and their incident edges.  

A {\it tree} is a nonempty simple, connected graph that contains no cycles. In a simple graph $\Gamma$, a {\it ray} is a one-way infinite simple path and a {\it double ray} or {\it line} is a two-way infinite simple path. For us the {\it ends} (sometimes called {\it vertex-ends} for non-locally-finite graphs) of $\Gamma$ are equivalence classes on the set of rays, where rays $R_1, R_2$ lie in the same end if and only if there exists a ray $R$ in $\Gamma$ containing infinitely many vertices of $R_i$ for $i=1,2$.
In a tree $T$, there is a unique shortest path between any two vertices $v$ and $w$, denoted $[v,w]$   (or $[v,w)$, for   example, if  we wish to exclude $w$). For an arc or edge $e$ in $T$ the graph $T \setminus \{e\}$ has two connected components; these are called the {\it half-trees} associated with $e$.

Actions of a group $G$ on a set $X$ are from the left, with $Gx$ denoting the orbit of $x \in X$ under the action of $G$. We denote the stabiliser of $x$ by $G_x$, and for a subset $Y \subseteq X$ we write $G_{\{Y\}}$ (resp.~$G_{(Y)}$) for the setwise (resp.~pointwise) stabiliser of $Y$ in $G$. The action is {\it transitive} on $X$ if $X = Gx$ for some $x \in X$. The group of all permutations of $X$ is denoted $\sym(X)$. Subgroups of $\sym(X)$ in which all orbits of point stabilisers are finite are called {\it subdegree-finite}.
A group $G \leq \sym(X)$ acts {\it freely} or {\it semiregularly} if the stabiliser $G_x$ of any $x \in X$ is trivial; if $G$ is transitive and semiregular we say it is {\it regular}. 
If $G$ is transitive, then it is {\it primitive} if and only if the only $G$-invariant equivalence relations on $X$ are the trivial relation (where equivalence classes are singletons) or the universal relation (where $X$ is an entire equivalence class).

There is a natural topology on $G$ that can be obtained from the action of $G$ on $X$, called the {\it permutation topology}, in which a neighbourhood basis of the identity is taken to be the pointwise stabilisers of finite subsets of $X$. If we think of $X$ as a discrete space with elements of $G$ as maps from $X$ to $X$, then the topology is equal to the topology of pointwise convergence and the compact-open topology. 
Permutational properties of $G$ have topological ramifications. For example, $G$ is totally disconnected if and only if the action on $X$ is faithful, and if $G \leq \sym(X)$ is closed and subdegree-finite then all stabilisers $G_x$ are compact and open, so $G$ is a totally disconnected and locally compact group (henceforth, \tdlc) with compact open stabilisers. See \cite{Moller:PermTDLC} for a thorough guide to this topology. 
  Topological statements concerning $\sym(X)$ will always appertain to the permutation topology.  
Note that a topological group is {\it compactly generated} if there is a compact subset that abstractly generates the group.

A {\it graph homomorphism} $\theta: \Gamma \rightarrow \Gamma'$ is a pair of maps $\theta_V: V\Gamma \rightarrow V\Gamma'$ and $\theta_A: A\Gamma \rightarrow A\Gamma'$ that respect origin vertices and edge reversal; if $\theta_V$ and $\theta_A$ are both bijections we say $\theta$ is an {\it isomorphism}. 
A group $G$ acting on $\Gamma$ gives rise to a {\it quotient graph} $G \backslash \Gamma$ whose vertex (resp.~ arc) set is the set of $G$-orbits on $V$ (resp.~on $A$), and for an arc $Ga$ in $G \backslash \Gamma$ we have $o(Ga) = Go(a)$ and $t(Ga) = Gt(a)$ and $\overline{Ga} = G\overline{a}$.
The group of automorphisms of $\Gamma$ is denoted $\Aut(\Gamma)$.  When $\Gamma$ is a simple graph, $\aut(\Gamma)$ acts faithfully on $V$ as those elements in $\sym(V)$ that respect the arc relation in $V\times V$, and in this case we identify $\Aut(\Gamma)$ with the corresponding subgroup of $\Sym(V)$.

For a tree $T$ (which recall is always simple) and a line $L$ in $T$, a {\it translation} of $L$ is an orientation-preserving automorphism of $L$ that does not fix any point on the line.
If $B$ is a subtree of $T$ and $G \leq \aut(T)$, we say that $G$ leaves $B$ {\it invariant} if $G$ fixes setwise the vertices of $B$.
Throughout, countable means finite or countably infinite. We say that an action on a tree $T$ is {\it geometrically dense} if the action does not leave invariant any nonempty proper subtree of $T$ and does not fix any end of $T$.

% ================================
% ================================
% ================================
% ================================
% ================================
% ================================
\section{Groups acting on trees}

%
%%
%%%
%%%%
%%%%%
%%%%%%
%%%%%%%
%%%%%%%%
%%%%%%%%%
%%%%%%%%%%
%%%%%%%%%%
%%%%%%%%%%
%%%%%%%%%%
\subsection{An overview of Bass--Serre Theory}
\label{Sec:BassSerre}

Traditionally, Bass--Serre Theory concerns groups acting on trees {\em without inversion}, meaning that no element of the group acts as an edge inversion. This is not a significant restriction, since given an action with inversion one can subdivide the edges of the tree and thus obtain an action without inversion. We abide by this tradition here and restrict our attention to inversion-free actions.   Let $T$ be a tree and let $G$ act on $T$ without inversion.  
We closely follow Section 5 of Serre's book \cite{Serre:trees}, so readers seeking a more complete description of the theory can consult this source. Our aim in this section is for readers to see that a local action diagram and its corresponding universal group are thoroughly dissimilar to a graph of groups and its universal cover, but nevertheless the beautiful correspondence in Bass--Serre Theory between the action and its description as a graph of groups is mirrored in our theory of local action diagrams. In Serre's work, a loop always admits an automorphism of order two which changes its orientation. We temporarily adopt this convention for Section~\ref{Sec:BassSerre}.

The combinatorial structure at the heart of Bass--Serre Theory is called a {\em graph of groups}. This is a connected nonempty graph $\Gamma$, together with some groups that will be associated with the vertices and edges of $\Gamma$; this association can be thought of as the vertices and edges of $\Gamma$ being coloured with a pallet of colours comprised of these groups. More precisely, for each vertex $v \in V\Gamma$ we have a group $\G_v$ (these are called the {\it vertex groups}) and for each arc $a \in A\Gamma$ we have a group $\G_a$ satisfying $\G_{\overline{a}} = \G_a$ (these are called the {\it edge groups}); we also have a monomorphism $\G_a \rightarrow \G_{t(a)}$ allowing us to view any edge group as a subgroup of the vertex groups of the vertices that comprise the edge. Let us denote this graph of groups as $(\Gamma, (\G_v), (\G_a))$.

This concise combinatorial structure admits two natural universal objects. The first is a tree called the {\em universal cover} of the graph of groups $(\Gamma, (\G_v), (\G_a))$; the second is the {\em fundamental group} of the graph of groups. Before formally defining these, let us first explore their significance. 
The first significant part of the Fundamental Theorem of Bass--Serre Theory (see Theorem~\ref{thm:FundThmBassSerre} below) is essentially a decomposition of the action $(T,G)$ in the language of graphs of groups:

\begin{quote} ($\circledast$) {\it There is a graph of groups associated to $(T,G)$, and $G$ can be identified with the fundamental group of this graph of groups.}
\end{quote}

The second significant part of the theorem is essentially a method of constructing inversion-free actions of groups on trees using graphs of groups:

\begin{quote} ($\circledcirc$) {\it Given a graph of groups $(\Gamma, (\G_v), (\G_a))$, its fundamental group acts on its universal cover (which is a tree) without inversion, and the graph of groups associated with this action   via $\circledast$ is   again $(\Gamma, (\G_v), (\G_a))$.}
\end{quote}

These two components of the theorem are the foundation of Bass--Serre Theory. For our purposes they are significant for two reasons. The constructive part ($\circledcirc$) is not usable in general for constructing new examples of nondiscrete actions of groups on trees because of what might be thought of as a `chicken or the egg' dilemma (see Section~\ref{subsection:bass_serre_limitations}). This limitation can be overcome via the local-to-global theory of groups acting on trees, and this local-to-global approach provided the motivation for our theory of local action diagrams. The second significance is that the rich correspondence between the action $(T,G)$ and the combinatorial object (i.e. the graph of groups) is mirrored in our theory of local action diagrams, in which the combinatorial object is a local action diagram.

Let us now formally define the objects in Bass--Serre Theory, and then give a precise statement of the theorem.

\subsubsection[The associated graph of groups]{The associated graph of groups (see \cite[\S 5.4]{Serre:trees})}
\label{sec:Assoc_graph}
Suppose $G$ acts on a tree $T$ without inversion. We first describe its {\it associated graph of groups}. The underlying graph of the graph of groups is the quotient graph $\Gamma := G \backslash T$. We next use $\Gamma$ to choose a subtree of $T$ in a coherent way so that we can take the vertex and arc stabilisers of this subtree to be the vertex and edge groups of the graph of groups.

Choose an orientation $E^+$ of $\Gamma$, and for each $a \in A\Gamma$ set $e(a) = 0$ whenever $a \in E^+$ and set $e(a) = 1$ otherwise.
Consider the subgraphs of $\Gamma$ that are trees; these form an ordered set (ordered by inclusion) and by Zorn's Lemma this set has a maximal element $M$, called a {\it maximal tree of $\Gamma$}, which is easily seen to contain every vertex of $\Gamma$ (see \cite[\S2 Proposition 11]{Serre:trees} for example). One can lift $M$ to a subtree $T'$ of $T$; that is, $T'$ is isomorphic to $M$ via the natural map $\pi$ which takes each vertex $v \in VT'$ (resp.~arc $a \in AT'$) to the vertex $Gv$ (resp.~arc $Ga$) in $M \subseteq G \backslash T$. 

Next, we construct a map $\varphi$ that takes arcs in $A\Gamma$ to arcs in $AT$ with the property that $\varphi(\overline{a}) = \overline{\varphi(a)}$ for all $a \in A\Gamma$. The map will also be defined on the vertices of $M$. On $M$ take $\varphi$ to be $\pi\inv: M \rightarrow T'$. For $a \in E^+ \setminus AM$ there is an arc $b \in AM$ such that $o(a) = o(b)$, and we set $\varphi(a) = \varphi(b)$. In this way we have $o(\varphi(a)) = \varphi(o(a))$. Now $t(\varphi(a))$ and $\varphi(t(a))$ have the same image under $\pi$ since both project to $t(b) \in V\Gamma$. In particular, $t(\varphi(a))$ and $\varphi(t(a))$ lie in the same $G$ orbit, so we can choose $\gamma_a \in G$ such that $t(\varphi(a)) = \gamma_a \varphi(t(a))$. From this we obtain elements $\gamma_a \in G$ for all arcs $a \in A\Gamma$ by specifying that $\gamma_{\overline{a}} = \gamma_a^{-1}$ and $\gamma_a$ is the identity whenever $a \in AM$.

We can now define the vertex and edge groups of our associated graph of groups. For each vertex $v$ of $\Gamma$ 
  recall that $v$ is a vertex in $M$, so $\varphi(v)$ is defined and is equal to a vertex in $T$;  
take the vertex group $\G_v$ to be the vertex stabiliser $G_{\varphi(v)} \leq G$. Similarly, for an arc $a$ of $\Gamma$ recall that $\varphi(a)$ is an arc in $T$; take the edge group $\G_a$ to be the arc stabiliser $G_{\varphi(a)} \leq G$, with $\G_{\overline{a}} = \G_a$. Finally, we define the monomorphism $\G_a \rightarrow \G_{t(a)}$ as $g \mapsto \gamma_a^{e(a)-1} g \gamma_a^{1-e(a)}$.

\begin{example}\label{Ex:AssocGraphOfGroups}    
Let $T$ be the biregular tree $T_{m,n}$ for finite distinct $m, n > 2$. Notice that the action $(T, \aut(T))$ has one edge orbit and two vertex orbits, and these vertex orbits correspond to the natural bipartition of $T_{m,n}$ into vertices with valency $m$, and vertices with valency $n$.  Thus, the graph of groups associated with $(T, \aut(T))$ is a pair of vertices connected by a single edge. The vertex groups are the stabilisers $\aut(T)_v$ and $\aut(T)_w$ where $v,w$ are adjacent vertices in $VT$ with $v$ having valency $m$ in $T$ and $w$ having valency $n$, and the edge group is $\aut(T)_{(v,w)}$.
\end{example}  

\subsubsection[The fundamental group of a graph of groups]{The fundamental group of a graph of groups (see \cite[\S 5.1]{Serre:trees})} \label{Sec:BassSerreFundGp}

Suppose we have a graph of groups $(\Gamma, (\G_v), (\G_a))$. For each 
arc $a \in A\Gamma$ we have a monomorphism $\G_a \rightarrow \G_{t(a)}$, and we denote the image of any $h \in \G_a$ under this monomorphism by $h^a$.
  We need extra generators for the fundamental group, one for each arc in $\Gamma$, so for all $a \in A\Gamma$  
we find new letters $g_a$ not contained in any of the vertex or edge groups. Our generating set $\Omega$ is   then the   union of $\{g_a : a \in A\Gamma\}$ and the vertex groups $\bigcup_{v \in V\Gamma} \G_v$. Let $M$ be a maximal tree of $\Gamma$. Then we define the {\it fundamental group $\Pi_1((\Gamma, (\G_v), (\G_a))$}, abbreviated to $\Pi_1$, to be
\[\langle \Omega : g_{\overline{a}} = g_a^{-1}, \ g_a h^a g_a^{-1} = h^{\overline{a}} \ (\forall a \in A\Gamma, h \in \G_a), \ g_b = 1 \ (\forall b \in AM) \, \rangle.\]

One can show (see \cite[\S 5 Proposition 20]{Serre:trees}) that the definition of $\Pi_1$ is independent of the choice for $M$.

\begin{example} \label{Ex:FundGp}
If $(\Gamma, (\G_v), (\G_a))$ is the graph of groups in Example~\ref{Ex:AssocGraphOfGroups}, then the definition above gives $\Pi_1$ to be the amalgamated free product $\aut(T)_v \ast_{\aut(T)_{(v,w)}} \aut(T)_w$.
\end{example}

\subsubsection[Universal covering of a graph of groups]{Universal covering of a graph of groups (see \cite[\S 5.3]{Serre:trees})}
\label{Sec:UnivCoveringBassSerre}

Now suppose we are given a graph of groups $(\Gamma, (\G_v), (\G_a))$, a maximal tree $M$ of $\Gamma$ and an orientation $E^+$ of $\Gamma$. Let $\Pi_1$ be the fundamental group of $(\Gamma, (\G_v), (\G_a))$.
Define the map $e : A\Gamma \rightarrow \{0,1\}$ as in Section~\ref{sec:Assoc_graph}, with $e(a) = 0$ whenever $a \in E^+$ and $e(a) = 1$ otherwise. For each arc $a \in \Gamma$, let $\hat{a} \in \{a, \overline{a}\}$ be such that $\hat{a} \in E^+$, and let $\G_a^a$ be the image (in $\G_{t(a)}$) of the edge group $\G_a$ under the monomorphism $h \mapsto h^a$. Notice that $\hat{a} = \hat{\overline{a}}$.

The {\it universal cover} of $(\Gamma, (\G_v), (\G_a))$ is a graph $\tilde{T}$ whose vertex set is the disjoint union of left cosets,
\[V\tilde{T} := \bigsqcup_{v \in V\Gamma}   \Pi_1   / \G_v,\]
and whose arc set is the disjoint union,
\[\bigsqcup_{a \in A\Gamma}   \Pi_1   / \G_{\hat{\overline{a}}}^{\hat{\overline{a}}}.\]
Arc inversion in $\tilde{T}$, and the maps $o$ and $t$, are defined as follows. For each arc $a \in A\Gamma$, let $\tilde{a}$ denote the trivial coset in $  \Pi_1   / \G_{\hat{\overline{a}}}^{\hat{\overline{a}}}$ corresponding to $\G_{\hat{\overline{a}}}^{\hat{\overline{a}}}$, and
for each vertex $v \in V\Gamma$, let $\tilde{v}$ denote the trivial coset in $  \Pi_1   / \G_{v}$ corresponding to $\G_{v}$.
Any $\tilde{T}$-arc lies in $  \Pi_1   / \G_{\hat{\overline{a}}}^{\hat{\overline{a}}}$ for some arc $a \in A\Gamma$.
Each $\tilde{T}$-arc in $  \Pi_1   / \G_{\hat{\overline{a}}}^{\hat{\overline{a}}}$ is of the form $g \tilde{a}$, for some $g \in   \Pi_1$.   Recalling our elements $g_a \in   \Pi_1$   from Section~\ref{Sec:BassSerreFundGp} and setting $v_o := o(a)$ and $v_t := t(a)$, we then take
\[\overline{g\tilde{a}} = g \tilde{\overline{a}}, \quad o(g\tilde{a}) = gg_a^{-e(a)} \widetilde{v_o}, \quad t(g\tilde{a}) = gg_a^{1-e(a)} \widetilde{v_t}.\]

One can then (see \cite[\S 5.3]{Serre:trees}) check that: $\tilde{T}$ is a tree and   $\Pi_1$   acts on the tree $\tilde{T}$ via left multiplication, and moreover this action is as automorphisms of $\tilde{T}$. Under this action the quotient graph $  \Pi_1   \backslash \tilde{T}$ is $\Gamma$, and for all vertices $v \in V\Gamma$ (resp.~arcs $a \in A\Gamma$) we have that the stabiliser $\Pi_{1 \tilde{v}}$  (resp.~$\Pi_{1 \tilde{a}}$) is equal to $\G_v$ (resp.~$\G_{\hat{\overline{a}}}^{\hat{\overline{a}}} \cong \G_a$). For any arc $a$ in the maximal tree $M$ of $\Gamma$ we have $g_a = 1$ and thus $\widetilde{o(y)} = o(\tilde{y})$ and $\widetilde{t(y)} = t(\tilde{y})$, so we have a lift of $M$ to a subtree $\tilde{T}'$ of the tree $\tilde{T}$ via $v \mapsto \tilde{v}$ and $a \mapsto \tilde{a}$ for $v \in VM$ and $a \in AM$. Thus, the associated graph of groups for $(\tilde{T},   \Pi_1)$   is $(\Gamma, (\G_v), (\G_a))$.

\begin{example} \label{Ex:UnivCover}   Let us continue Examples~\ref{Ex:AssocGraphOfGroups}--\ref{Ex:FundGp}, resuming their notation. For the tree $T = T_{m,n}$ and the action $(T, \aut(T))$, the graph of groups $(\Gamma, (\G_v), (\G_a))$ consists of two vertices connected by a single edge, with vertex groups $\aut(T)_v$ and $\aut(T)_w$, and edge group $\aut(T)_{v,w}$, where $v$ has valency $m$ in $T$ and $w$ has valency $n$. 
Recall that the fundamental group of this graph of groups is $\Pi_1=\aut(T)_v \ast_{\aut(T)_{(v,w)}} \aut(T)_w$.
Now $\aut(T)_v$ is transitive on the edges in $T$ that are incident to $v$, so the index of the edge group $\aut(T)_{(v,w)}$ in the vertex group $\aut(T)_v$ is $m$. Similarly the index of $\aut(T)_{(v,w)}$ in $\aut(T)_w$ is $n$.
The universal cover $\tilde{T}$ of $(\Gamma, (\G_v), (\G_a))$ is a tree, and by definition we have $V\tilde{T} = \left ( \Pi_1 / \aut(T)_v \right ) \sqcup \left ( \Pi_1  / \aut(T)_w \right )$ with the two sets in this disjoint union naturally bipartitioning $T$.  
\end{example}

\subsubsection[Fundamental Theorem of Bass--Serre Theory]{Fundamental Theorem of Bass--Serre Theory (see \cite[\S5.4]{Serre:trees})}
\label{Sec:FundThmBassSerre}

We are now able to formally state the fundamental theorem.

\begin{theorem}[{\cite[\S5]{Serre:trees}}] \label{thm:FundThmBassSerre}
Let $G$ be a group and $T$ be a tree.

{\normalfont ($\circledast$)} Suppose $G$ acts on $T$ without inversion. Let $(\Gamma, (\G_v), (\G_a))$ be its associated graph of groups, and choose a maximal tree $M$ of $\Gamma$. Let $  \Pi_1$   be the fundamental group of this graph of groups (with respect to $M$) and let $\tilde{T}$ be the universal covering (with respect to $M$). Then the map $\phi :   \Pi_1   \rightarrow G$ defined by the inclusion $\G_v \leq G$ (for $v \in V\Gamma$) and $\phi(g_a) = \gamma_a$ (for $a \in A\Gamma$) is an isomorphism of groups. The map $\psi : \tilde{T} \rightarrow T$ given, for all $h \in   \Pi_1$  , by $\psi(h \tilde{v}) = \phi(h) \varphi(v)$ (for all $v \in V\Gamma$) and $\psi(h\tilde{a}) = \phi(h) \varphi(a)$ (for all $a \in A\Gamma$) is an isomorphism of graphs. Moreover, the isomorphism $\psi$ is $\phi$-equivariant, and so the actions $(\tilde{T},   \Pi_1)$   and $(T,G)$ can be identified.

{\normalfont ($\circledcirc$)} Suppose $(\Gamma, (\G_v), (\G_a))$ is a graph of groups. Let   $\Pi_1$   be its fundamental group and $\tilde{T}$ its universal cover. Then $(\tilde{T},   \Pi_1)  $ is an inversion-free action on a tree whose associated graph of groups is $(\Gamma, (\G_v), (\G_a))$.
\end{theorem}

\begin{example}  
Let us continue Examples~\ref{Ex:AssocGraphOfGroups}--\ref{Ex:UnivCover}, resuming their notation. We have the tree $T = T_{m,n}$ and a graph of groups $(\Gamma, (\G_v), (\G_a))$ for the action $(T, \aut(T))$. The universal cover of $(\Gamma, (\G_v), (\G_a))$ is $\tilde{T}$ and its fundamental group is $\aut(T)_v \ast_{\aut(T)_{(v,w)}} \aut(T)_w$. By the Fundamental Theorem of Bass--Serre Theory, we have that the actions of $\aut(T)_v \ast_{\aut(T)_{(v,w)}} \aut(T)_w$ on $\tilde{T}$ and of $\aut(T)$ on $T$ are permutationally isomorphic.  
\end{example}

%
%%
%%%
%%%%
%%%%%
%%%%%%
%%%%%%%
%%%%%%%%
%%%%%%%%%
%%%%%%%%%%
%%%%%%%%%%
%%%%%%%%%%
%%%%%%%%%%
\subsection{Jacques Tits' Independence Property $\propP{}$ and simplicity}
\label{Section:tits}

We describe \cite[\S 4.2]{Tits70}, where   {\it Tits' Independence Property $\propP{}$} (also called {\it property $\propP{}$} or the {\it independence property})   is introduced.
Suppose $G \leq \aut(T),$   where $T$ is a tree.   If $C$ is a (finite or infinite) nonempty simple path in $T$, then for each $v \in VT$ there is a unique vertex $\pi_{C}(v)$ in $C$ which is closest to $v$. This gives a well-defined map $v \mapsto \pi_{C}(v)$. For each $w \in VC$, the set $\pi^{-1}_{C}(w)$ is the vertex set of a subtree of $T$. Each of these subtrees is invariant under the action of the pointwise stabiliser $G_{(C)}$ of $C$, and so we define $G^w_{(C)}$ to be the subgroup of $\sym(\pi^{-1}_{C}(w))$ induced by $G_{(C)}$. We therefore have homomorphisms $\varphi_w: G_{(C)} \rightarrow G_{(C)}^w$ for each $w \in VC$ from which we obtain the homomorphism,
\begin{equation*}
\label{equation:propP}
	\varphi: G_{(C)} \rightarrow \prod_{w \in VC} G_{(C)}^w.
\end{equation*}
Now $G$ has the   {\it independence property for $C$}   if the homomorphism $\varphi$ is an isomorphism,   and $G$ has {\it Tits' Independence Property $\propP{}$} if it has the independence property for $C$   for every possible simple path $C$.
Intuitively, $G$ has   the independence property for a path $C$ if $G_{(C)}$ can act independently on all of the subtrees `hanging' from $C$. 

\begin{theorem}[{\cite[Th\'{e}or\`{e}me 4.5]{Tits70}}] \label{theorem:Tits4_5}
Let $T$ be a tree and suppose $G \leq \aut(T)$ has property $\propP{}$. Let $G^+$ be the group generated by the pointwise stabilisers in $G$ of edges in $T$. If $G$ is geometrically dense, then $G^+ \unlhd G$ and $G^+$ is trivial or abstractly simple.
\end{theorem}

%
%%
%%%
%%%%
%%%%%
%%%%%%
%%%%%%%
%%%%%%%%
%%%%%%%%%
%%%%%%%%%%
%%%%%%%%%%
%%%%%%%%%%
%%%%%%%%%%
\section{Universal groups acting on trees}

%
%%
%%%
%%%%
%%%%%
%%%%%%
%%%%%%%
%%%%%%%%
%%%%%%%%%
%%%%%%%%%%
%%%%%%%%%%
%%%%%%%%%%
%%%%%%%%%%
\subsection{Limitations to Bass--Serre Theory}
\label{subsection:bass_serre_limitations}

In a graph of groups, we specify embeddings of edge groups into vertex groups. In doing so, the vertex stabiliser is specified {\it as an abstract group}, while the action on the $1$-ball is also specified. In many commonly encountered situations, this leads to complications that are intractable. We describe three of them here; they are related but different.

1. In the action on a tree resulting from a graph of groups, the vertex stabiliser can only be a quotient of the vertex group. So, we can't obtain `large' (e.g.~nondiscrete) vertex stabilisers this way, unless the vertex group we started with already had an interesting action on a tree. 

2. Suppose we are given a graph of groups and in addition we are told that the action is in some sense interesting. Even though we have been given the vertex stabiliser, how it acts on the tree on a large scale is unintuitive and often abstruse. Even basic questions, like whether or not the action is faithful, can be difficult to answer. Restricting to locally finite trees does not sufficiently reduce the complexity: we are still in a sense required to understand the dynamics of sequences of virtual isomorphisms as we travel along all possible walks through the graph of groups. This is typically an insurmountable problem when the vertex stabilisers are infinite, or when one is interested in understanding ever larger finite vertex stabilisers as in the Weiss conjecture.

3. The possible combinations of vertex groups are difficult to understand because they are not independent of one another: we need a common edge group for any pair of neighbouring vertices. If we wish to choose a collection of local actions independently of one another, there are very few ways to do this and all have significant limitations. For example, we could make all the vertex groups infinitely generated free groups, but then determining the action of the resulting vertex groups on the tree is hopelessly complicated, rendering the task of understanding the group's closure in $\aut(T)$ as a topological group futile.

Readers who are familiar with Bass--Serre Theory will no doubt recognise these limitations. For readers unfamiliar with Bass--Serre Theory, we set an (intractable) exercise that highlights some of these issues.

\begin{exercise*} 
Let $T$ be the biregular tree $T_{7,5}$. Using only Bass--Serre Theory,   attempt to   construct a subgroup $G \leq \aut(T)$ that has two vertex orbits on $T$, such that every 
vertex stabiliser $G_v$ is infinite and does not induce $S_7$ or $S_5$ on the set of neighbours of the vertex $v$.   The point of this exercise is not to complete the construction, it is to make an attempt and in doing so experience the aforementioned limitations 1--3.  
\end{exercise*}

As we shall see, the exercise has an easy solution using the theory of local action diagrams. However, even when we take this solution and (pretending for a moment that we do not know it is a solution) use Bass--Serre Theory to analyse its action, our point (2) above becomes apparent.

%
%%
%%%
%%%%
%%%%%
%%%%%%
%%%%%%%
%%%%%%%%
%%%%%%%%%
%%%%%%%%%%
%%%%%%%%%%
%%%%%%%%%%
%%%%%%%%%%
\subsection{  Burger--Mozes   groups}
\label{Section:BurgerMozes}

Here we largely follow Burger and Mozes' paper \cite[\S 3.2]{BurgerMozes} with one significant difference: we do not insist that the trees be locally finite. Let $d > 2$ be some finite or infinite cardinal and let $T$ be a regular tree of valency $d$. Fix some set $X$ such that $|X| = d$, and let $F \leq \sym(X)$. We say that $G \leq \aut(T)$ is {\it locally-$F$} if for all $v \in VT$ the permutation group induced on the set $B(v)$ of neighbours of $v$ by the vertex stabiliser $G_v$ is permutationally isomorphic to $F$.
A {\it legal colouring} is a map $\lcol: AT \rightarrow X$ that satisfies the following:
\begin{enumerate}%[label=(\roman*)]
    \item
        For each vertex $v \in VT$, the restriction 
        $\lcol \big|_{o\inv(v)}: o\inv(v) \rightarrow X$ is a bijection;
    \item
        For all arcs $a \in AT$ we have $\lcol(a) = \lcol(\overline{a})$.
\end{enumerate}
One can always construct a legal colouring of $T$.
For $g \in \aut(T)$ and $v \in VT$ we define the {\it $\lcol$-local action of $g$ at $v$} to be,
\begin{equation}\label{eq:LocalAction}
\sigma_{\lcol, v}(g):= \lcol \big|_{o\inv(gv)} g \lcol \big|^{-1}_{o\inv(v)}.
\end{equation}
Notice that $\sigma_{\lcol, v}(g) \in \sym(X)$. One might ask if we can constrain these bijections, so that they all lie in some common subgroup of $\sym(X)$. Such a restriction gives rise to the Burger--Mozes universal groups.

The {\it Burger--Mozes universal group of $F$} (with respect to $\lcol$) is the group,
\[\U_\lcol(F) := \{g \in \aut(T) : \sigma_{\lcol, v}(g) \in F \quad \forall v \in VT\}.\]
Two legal colourings $\lcol$ and $\lcol'$ give rise to universal groups $\U_\lcol(F)$ and $\U_{\lcol'}(F)$ that are conjugate in $\aut(T)$; for this reason we replace $\U_\lcol(F)$ with $\U(F)$ and speak of {\it the} Burger--Mozes universal group $\U(F)$ of $F$. Further properties of $\U(F)$ described in \cite{BurgerMozes}  are as follows. Note that some of these properties are expanded upon in Section~\ref{Section:BoxProduct}.
\begin{enumerate}
    \item
        $\U(F)$ is a vertex transitive and locally-$F$ subgroup of $\aut(T)$, and if $F$ has finite degree then $\U(F)$ is closed;
    \item
        $\U(F)$ enjoys Tits Independence Property $\propP{}$, and consequently by Theorem~\ref{theorem:Tits4_5} the subgroup $\U(F)^+$ is trivial or simple;
    \item
        If $F$ has finite degree then the index $|\U(F) : \U(F)^+|$ is finite if and only if $F \leq \sym(X)$ is transitive and generated by point stabilisers; when this happens $\U(F)^+ = \U(F) \cap (\aut(T))^+$ and therefore $|\U(F) : \U(F)^+| = 2$;
    \item
        when $F$ is transitive on $X$, the group $\U(F)$ has the following universal property: $\U(F)$ contains an $\aut(T)$-conjugate of every locally-$F$ subgroup of $\aut(T)$.
\end{enumerate}
Thus, in the situation when $F$ is transitive, we have a natural way to describe $\U(F)$: it is the largest locally-$F$ subgroup of $\aut(T)$.

Burger and Mozes used $U(F)$ to build towards a hoped-for classification of closed $2$-transitive actions on locally finite trees. Here we instead see them as a first step towards a local-to-global theory of groups acting on trees.
When $T$ is locally finite, $\U(F)$ is obviously a \tdlc subgroup of $\aut(T)$ with compact vertex stabilisers. In fact, the group can be \tdlc in more general situations; this was shown in \cite{SmithDuke}, where the Burger--Mozes universal groups are viewed as a special case of a more general construction called the {\it box product}. In this more general setting, various global topological properties are easy to characterise using local conditions (i.e.~conditions on $F$). We describe these in section \ref{Section:BoxProduct}.

%
%%
%%%
%%%%
%%%%%
%%%%%%
%%%%%%%
%%%%%%%%
%%%%%%%%%
%%%%%%%%%%
%%%%%%%%%%
%%%%%%%%%%
%%%%%%%%%%
\subsection{The box product construction}
\label{Section:BoxProduct}

Here we largely follow the second author's paper \cite{SmithDuke}, where the box product was introduced as a product for permutation groups; as a permutational product it is, in some sense, the dual of the wreath product in its product action. 
In the same paper, the box product was used to produce the first examples of $2^{\aleph_0}$ distinct isomorphism classes of compactly generated locally compact groups that are nondiscrete and topologically simple, answering a long standing open question. The box product arises naturally in the structure theory of subdegree-finite primitive permutation groups due to the second author (\cite{SmithPrimitive}).

Let $d_1, d_2 > 1$ be two finite or infinite cardinal numbers and let   $T = T_{d_1, d_2}$ be the   biregular tree with valencies $d_1, d_2$.
Fix disjoint sets $X_1, X_2$ such that $|X_1| = d_1$ and $|X_2|=d_2$. Thus there is a natural bipartition of $T$ as $VT = V_1 \sqcup V_2$ such that vertices in $V_1$ have valency $|X_1|$ and vertices in $V_2$ have valency $|X_2|$.
Let $F_i \leq \sym(X_i)$, for $i = 1,2$ with at least one of $F_1, F_2$ being nontrivial. We say $G \leq \aut(T)$ is {\it locally-($F_1, F_2$)} if $G$ preserves setwise the parts $V_1$ and $V_2$, and for all $v \in VT$ 
the permutation group induced on $B(v)$ by the vertex stabiliser $G_v$ is permutationally isomorphic to $F_i$ whenever $v \in V_i$, for $i=1,2$.

In this new context we define a {\it legal colouring} to be a map $\lcol: AT \rightarrow X$ that satisfies the following:
\begin{enumerate}
    \item
        For each vertex $v \in V_i$, the restriction 
        $\lcol \big|_{o\inv(v)}: o\inv(v) \rightarrow X_i$ is a bijection, for $i=1,2$;
    \item
        For all $v \in VT$ the restriction $\lcol \big|_{t\inv(v)}$ is constant.
\end{enumerate}
One can always construct a legal colouring of $T$. Using the same $\lcol$-local action defined in Equation~\ref{eq:LocalAction}, we define the {\it topological box product of $F_1$ and $F_2$} to be the group
\[\U_\lcol(F_1, F_2) := \left \{g \in \aut(T)_{\{V_1\}} : \sigma_{\lcol, v}(g) \in   F_i \quad \forall v \in V_i, \text{ for } i=1,2   \right \}.\]
The {\it permutational box product of $F_1$ and $F_2$} is the subgroup of $\sym(V_2)$ induced by the action of $\U_\lcol(F_1, F_2)$ on $V_2$, and is denoted $F_1 \boxtimes_{\lcol} F_2$.

As with the Burger--Mozes groups, two legal colourings $\lcol$ and $\lcol'$ give rise to groups $\U(F_1, F_2)_\lcol$ and $\U(F_1, F_2)_{\lcol'}$ that are conjugate in $\aut(T)$ and so we speak of {\it the} box product and write $\U(F_1, F_2)$ and $F_1 \boxtimes F_2$.

The Burger--Mozes groups arise as special cases of the box product construction: for any permutation group $F \leq \sym(X)$ where $|X| \geq 3$, the Burger--Mozes group $\U(F)$ is permutationally isomorphic to $S_2 \boxtimes F$, where $S_2$ here denotes the symmetric group of degree $2$. To see why, write $d := |X|$ and note that $\U(F)$ as a group of automorphisms of the $d$-regular tree $T_d$ induces a faithful action on $T_{2,d}$ because $T_{2,d}$ is the barycentric subdivision of $T_d$.
Now $\U(S_2, F)$ also acts as a group of automorphisms on $T_{2,d}$ and one can easily verify that the actions of $\U(F)$ and $\U(S_2, F)$ induced on the set of $d$-valent vertices of $T_{2,d}$ is permutationally isomorphic. From this we have that $\U(F)$ is topologically isomorphic to $S_2 \boxtimes F$, $\U(F, S_2)$, $F \boxtimes S_2$, and $\U(S_2, F)$.\\

As a permutational product, the box product has the following striking similarity to the wreath product in its product action. Recall that $F_1 \Wr F_2$ (in its product action) is a primitive permutation group if and only if $F_1$ is primitive and nonregular and $F_2$ is transitive and finite. Special cases of this fact were, according to Peter M.~Neumann, first proved by W.~A.~Manning in the early 20th Century. For over a century, no other product of permutation groups was known to preserve primitivity in this kind of generality. Compare this with the box product: $F_1 \boxtimes F_2$ is a primitive permutation group if and only if $F_1$ is primitive and nonregular and $F_2$ is transitive.

With the benefit of hindsight, we can see hints of this local and global primitivity equivalence for the box product construction going back to the 1970s. Indeed, groups of automorphisms of simple graphs with connectivity one (see \cite[\S 7.2]{ReidSmith}) can be realised as faithful inversion-free actions on trees, and in 1977 H.~A.~Jung and M.~E.~Watkins in \cite[Theorem 4.2]{JungWatkins} proved that the automorphism group of a simple graph $\Gamma$ with connectivity one is (i) vertex primitive if and only if (ii) all its lobes are vertex primitive, have at least 3 vertices and are pairwise isomorphic. Now (i) implies that the induced faithful tree action is primitive on one part $P$ of the bipartition of the tree, and (ii) implies that the local actions at vertices in $T \setminus P$ are primitive. 
Arguments by R\"{o}gnvaldur G. M\"{o}ller in the 1994 paper \cite{Moller:Primitive} (which used Warren Dicks and M.~J.~Dunwoody's powerful theory of structure trees from \cite{DicksDunwoody}) can be used to show all infinite, subdegree-finite primitive permutation groups with more than one end have a locally finite orbital digraph with connectivity one.
In 2010 Jung and Watkins' result was generalised to directed graphs by the second author in \cite{SmithDigraph}, and in this context the primitivity condition in (ii) becomes primitive but not regular. These observations inspired constructions in the second author's 2010 paper \cite{SmithSubdegrees} which, again with hindsight, can be seen as precursors to the permutational box product construction. As noted in \cite{SmithDuke}, the box product $F_1 \boxtimes F_2$ can be constructed for closed groups $F_1, F_2$ using refinements of the arguments in \cite{SmithSubdegrees} (the arguments are based on countable relational structures).

We also see hints coming from the world of groups acting on   locally finite   trees: obviously, in Burger and Mozes' 2000 paper \cite{BurgerMozes}, but also in Pierre-Emmanuel Caprace and Tom De Medts' 2011 result \cite[Theorem 3.9]{CapraceDeMedts}, which states the following. 
  Let $T$ be a locally finite tree and suppose $G \leq \aut(T)$ is nondiscrete, noncompact, compactly generated, closed, \tdlc and topologically simple with Tits' Independence Property $\propP{}$. Then the following conditions are equivalent: 
(i) every proper open subgroup of $G$ is compact; and (ii) $G$ splits as an amalgamated free product $G \cong A \ast_C B$, where $A$ and $B$ are maximal compact open subgroups, $C = A \cap B$ and the $A$-action on $A/C$ (resp.~the $B$-action on $B/C$) is primitive and noncyclic.  

For a group $G$   satisfying Caprace and De Medts' result,   condition (i) implies $G$ has a natural permutation representation that is primitive, and (ii) implies that the action of $G$ on its Bass--Serre tree is primitive on {\it both} parts of its bipartition and moreover the local action at vertices in both parts of the bipartition is primitive but not regular.\\

The following are further properties of the box product construction. Since $\U(F)$ is topologically isomorphic to $\U(F, S_2)$ and permutationally isomorphic to $S_2\boxtimes F$, many of these properties expand the properties of Burger--Mozes groups given in Section~\ref{Section:BurgerMozes}.
\begin{enumerate}
    \item
        $\U(F_1, F_2)$ is a locally-$(F_1, F_2)$ subgroup of $\aut(T)$.
    \item
        Any subset $Y \subseteq X_1 \cup X_2$ is an orbit of $F_1$ or $F_2$ if and only if $t(\lcol\inv (Y))$ is an orbit of $\U_\lcol(F_1, F_2)$. In particular, $F_1 \boxtimes F_2$ is transitive if and only if $F_1$ is transitive, and $\U(F_1, F_2)$ has precisely two vertex-orbits ($V_1$ and $V_2$) if and only if $F_1$ and $F_2$ are transitive.
    \item
        If $F_i$ is a closed subgroup of $\sym(X_i)$ for $i=1,2$ then $\U(F_1, F_2)$ is a closed subgroup of $\aut(T)$.
    \item
        $\U(F_1, F_2)$ enjoys Tits Independence Property $\propP{}$.
    \item \label{item:BoxSimple}
        If $F_1, F_2$ are generated by point stabilisers, then $\U(F_1, F_2)$ is simple if and only if $F_1$ or $F_2$ is transitive.
    \item
        When $F_1, F_2$ are transitive, the group $\U(F_1, F_2)$ has the following universal property: $\U(F_1, F_2)$ contains an $\aut(T)$-conjugate of every locally-$(F_1, F_2)$ subgroup of $\aut(T)$.
    \item
        If $F_1, F_2$ are closed, then $\U(F_1, F_2)$ is locally compact (and hence \tdlc) if and only if all point stabilisers in $F_1$ and $F_2$ are compact. Moreover, for all $v \in V_2$ the stabiliser $\U(F_1, F_2)_v$ is compact if and only if $F_2$ is compact and every point stabiliser in $F_1$ is compact.
    \item
        If $F_1, F_2$ are closed with compact point stabilisers, then all point stabilisers in $\U(F_1, F_2)$ are compactly generated if and only if $F_1$ and $F_2$ are compactly generated. Moreover, if $F_1$ and $F_2$ are compactly generated with finitely many orbits, and at least one of the groups is transitive, then $\U(F_1, F_2)$ is compactly generated.
    \item
        $\U(F_1, F_2)$ is discrete if and only if $F_1$ and $F_2$ are semiregular.
\end{enumerate}
Thus, in the situation when $F_1, F_2$ are transitive, we again have a natural way to describe $\U(F_1, F_2)$: it is the largest locally-$(F_1, F_2)$ subgroup of $\aut(T)$. Furthermore, we can under mild {\it local} conditions (that is, conditions on $F_1, F_2$) ensure that $\U(F_1, F_2)$ is nondiscrete, compactly generated, \tdlc and abstractly simple. An important point here is that one can use discrete local groups $F_1, F_2$ with certain topological properties, and obtain a nondiscrete topological group $\U(F_1, F_2)$ that inherits those local topological properties (except of course, discreteness) as global topological properties, as well as being guaranteed other desirable global properties like Tits' Independence Property $\propP{}$ and nondiscreteness.

Constructing $2^{\aleph_0}$ distinct isomorphism classes of compactly generated locally compact groups
that are nondiscrete and topologically simple is now relatively straightforward. For example, for a large enough prime $p$ (e.g.~$p> 10^{75}$) let $\{\mathcal{O}_i\}_{i \in I}$ be a set of representatives from the $2^{\aleph_0}$ isomorphism classes of   A.~Yu.~Ol'shanski{\u\i}'s   $p$-Tarski Monsters (\cite{olshanski}).
Each $\mathcal{O}_i$ is infinite and simple and any nontrivial proper subgroup $H_i < \mathcal{O}_i$ has finite order $p$. The groups $\mathcal{O}_i$ can therefore be viewed as (faithful) groups of permutations of $\mathcal{O}_i / H_i$. The groups $\U_i := \U(\mathcal{O}_i, S_3)$ are nondiscrete, compactly generated, \tdlc and simple.
It can then be shown (with a little work) that the $\U_i$ are pairwise nonisomorphic.

%
%%
%%%
%%%%
%%%%%
%%%%%%
%%%%%%%
%%%%%%%%
%%%%%%%%%
%%%%%%%%%%
%%%%%%%%%%
%%%%%%%%%%
%%%%%%%%%%
\subsection{Further generalisations}
\label{Sec:FurtherGens}

In this section we survey a selection of generalisations of the Burger--Mozes groups, the box product construction, and Tits' Independence Property $\propP{}$. These generalisations play no part in the theory of local action diagrams, but they suggest natural generalisations to the theory that we will revisit in Section~\ref{Section:FutureWork}. \\

Adrien Le Boudec in \cite{LeBoudec} considered universal groups of the locally finite $d$-valent tree $T_d$, for $d \geq 3$, that have local action $F \leq S_d$ at all but finitely many vertices. (Earlier, the special case where $F$ is the alternating group of degree $d$ was examined in \cite{BCGM}.)
More precisely, for a Burger--Mozes legal colouring $\lcol$ of $T_d$ and a group $F \leq S_d$, the {\it Le Boudec group} $G_{\lcol}(F)$ is the group,
\[\{g \in \aut(T_d) : \sigma_{\lcol, v}(g) \in F \quad \text{for all but finitely many } v \in VT_d\}.\]
 For a given $g \in G_{\lcol}(F)$, the finitely many vertices $v$ for which $\sigma_{\lcol, v}(g) \not \in F$ are called the {\it singularities of $g$}. Of course we have $\U_\lcol(F) \leq G_\lcol(F) \leq \aut(T_d)$. There is a topology on $G_\lcol(F)$ 
such that the inclusion map of $\U_\lcol(F)$ into $G_\lcol(F)$ is continuous and open (in general this is not the topology inherited from $\aut(T_d)$). 
Under this topology $G_\lcol(F)$ is a \tdlc group, and $G_\lcol(F)$ is discrete if and only if $F$ is a semiregular permutation group.
 
 For permutation groups $F \leq F' \leq S_d$, Le Boudec defines a group $G_\lcol(F, F')$ now called the {\it restricted universal group}, with $G_\lcol(F, F') := G_\lcol(F) \cap \U_\lcol(F')$ under the topology induced from $G_\lcol(F)$.
Note that, while the groups $G_\lcol(F)$ and $G_\lcol(F, F')$ are subgroups of $\aut(T_d)$, in general they are not closed in $\aut(T_d)$.

The Le Boudec groups and restricted universal groups are an important source of examples of compactly generated locally compact simple groups that do not admit lattices\footnote{A {\it lattice} $\Lambda$ in a locally compact group $G$ is a discrete subgroup such that the quotient space $G/\Lambda$ admits a $G$-invariant finite measure; the lattice is {\it uniform} if the quotient space $G/\Lambda$ is compact.} Prior to Le Boudec's work, the only example of a lattice-free compactly generated locally compact simple group was Neretin's group $\AAut(T_d)$ of almost automorphisms (see \cite{BCGM}). In particular Le Boudec uses these constructions to show the very nice result that there exist compactly generated abstractly simple \tdlc groups $H\leq G$ with $H$ is cocompact in $G$ such that $G$ contains lattices but $H$ does not.\\

Waltraud Lederle in \cite{Lederle} simultaneously generalises Neretin's group $\AAut(T_d)$ and the Burger--Mozes groups for the locally finite tree $T_d$. Some of the terms defined in \cite{Lederle} have subsequently become known by different names, and we follow this more recent naming convention here.

Let us first define Neretin's group $\AAut(T_d)$. Finite subtrees of $T_d$ in which every vertex is either a leaf or of valency $d$ are called {\it complete}. An {\it almost automorphism} of $T_d$ is an isomorphism $g$ of rooted forests $g : T_d \setminus C_0 \rightarrow T_d\setminus C_1$ where $C_0, C_1$ are complete finite subtrees of $T_d$. If 
$g : T_d \setminus C_0 \rightarrow T_d\setminus C_1$ and $g' : T_d \setminus C'_0 \rightarrow T_d\setminus C'_1$ are almost automorphisms, we say they are equivalent if and only if they agree on $T_d \setminus C$ for some complete finite subtree $C$ satisfying $C_0 \cup C'_0 \subseteq C$. Equivalent almost automorphisms induce the same homeomorphism of the set $\partial T_d$ of ends of $T_d$; in other words they give rise to the same {\it spheromorphism} of $\partial T_d$. The spheromorphisms are thus the equivalence classes of almost automorphisms under this equivalence relation.
To multiply two spheromorphisms $[g]$ and $[g']$, choose representatives $\overline{g} \in [g]$ and $\overline{g}' \in [g']$ that are defined on the same forest $T_d \setminus C$ (where $C$ is a finite complete subtree), then take $[g][g']$ to be the equivalence class containing $\overline{g} \,\overline{g}'$. The set of spheromorphisms under this multiplication form a group $\AAut(T_d)$ called {\it Neretin's group}. It was shown in \cite{BCGM} that Neretin's group 
can be given a group topology by taking the conjugates of vertex stabilisers in $\Aut(T_d) \leq \AAut(T_d)$ as a subbasis of neighbourhoods of the identity (this makes sense because $\AAut(T_d)$ commensurates each stabiliser $\Aut(T_d)_v$). Under this topology, $\AAut(T_d)$ is a compactly generated \tdlc group that contains $\aut(T_d)$ as an open subgroup. Moreover, $\AAut(T_d)$ is simple and contains no lattice (see \cite[Theorem 1]{BCGM}).

Given $G \leq \aut(T_d)$, Waltraud Lederle's group $\mathcal{F}(G)$ consists of elements of $\AAut(T_d)$ that intuitively `locally look like' $G$. When $G$ is closed and has Tits' Independence Property $\propP{}$ (that is, when $G$ is $\propP{}$-closed), there is a unique group topology on $\mathcal{F}(G)$ such that the inclusion of $G$ into $\mathcal{F}(G)$ is continuous and open. Since $G$ is \tdlc the group $\mathcal{F}(G)$ is also a \tdlc group. The focus of \cite{Lederle} is when $G$ is taken to be a Burger--Mozes group $\U(F)$.

To define $\mathcal{F}(G)$ precisely, we first define a {\it $G$-almost automorphism}   of $T_d$   to be an almost automorphism $h : T_d \setminus C \rightarrow T_d \setminus C'$ such that for every component $\mathcal{T}$ of the forest $T_d \setminus C$ there is some $g \in G$ such that $g$ and $h$ agree on $\mathcal{T}$. The elements of $\mathcal{F}(G)$ are the equivalence classes of all $G$-almost automorphisms and one can easily see that $\mathcal{F}(G) \leq \AAut(T_d)$.

In \cite{Lederle} it is shown that if $F \leq S_d$ then $\mathcal{F}(\U(F))$ is a compactly generated \tdlc group and its commutator subgroup is open, abstractly simple and has finite index in $\mathcal{F}(\U(F))$. Moreover, if $F$ is a Young subgroup with   strictly fewer   than $d$ orbits, then this commutator subgroup is a nondiscrete compactly generated simple group without lattices.
\\

Jens Bossaert and Tom De Medts in \cite{JensTom} generalise the box product construction by defining universal groups over right-angled buildings where the local actions are again arbitrary permutation groups. As with the box product, there is no requirement for these local actions to have finite degree. Bossaert and De Medts characterise when their universal groups are locally compact, abstractly simple or have primitive action on the residues of the building.\\

An important generalisation of Tits' Independence Property $\propP{}$ is described by Christopher Banks, Murray Elder and George A.~Willis in \cite{BanksElderWillis}, where it is called {\it Property $\propP{k}$.} This definition is made only for locally finite trees. For groups $H \leq \aut(T)$ satisfying Property $\propP{k}$ there is a natural subgroup denoted $H^{+_k}$ that is trivial or simple whenever $H$ neither fixes an end nor setwise stabilises a proper subtree of $T$ (see \cite[Theorem 7.3]{BanksElderWillis}).
In the same paper they define something they call the $k$-closure of a group. In a different context this notion was independently described by Sam Shepherd in \cite{SamShep}. The term $k$-closure has a different and well-established meaning in permutation group theory due to Wielandt, and so we do not use the term; we instead call this notion the {\it $\propP{k}$-closure}   and define it in the following paragraph.   In both \cite{BanksElderWillis} and \cite{SamShep} the notion of $\propP{k}$-closure pertains only to locally finite trees where issues surrounding the closure or non-closure of local actions do not occur, because they are always closed. In fact, the definition of $\propP{k}$-closure can be adjusted slightly so as to work in a natural way for trees that are not locally finite, and this latter notion was introduced in our paper \cite{ReidSmith}. Here we follow \cite{ReidSmith}, and note that for locally finite trees the definitions coincide.

Let $T$ be a tree that may or may not be locally finite, and let $G \leq \aut(T)$. For $k \in \N$, the {\it $\propP{k}$-closure} of $G$, denoted $G^{\propP{k}}$, is the group consisting of all $g \in \Aut(T)$ such that for all $v \in VT$, and all finite subsets $\Phi \subseteq B_v(k)$, there exists $g_\Phi \in G$ such that $gw = g_\Phi w$ for every vertex $w \in \Phi$.  We say that the group $G$ is {\it $\propP{k}$-closed} if $G$ is equal to its $\propP{k}$-closure. Intuitively, we can think of the $k$-local action of $G$ as being its action on $k$-balls $B_v(k)$, and 
the $\propP{k}$-closure of $G$ is then the largest subgroup of $\aut(T)$ 
whose $k$-local action is equal to the closure of the $k$-local action of $G$. Note that all of the groups described thus far in this subsection that arise as closed groups of automorphisms of trees 
  are $\propP{k}$-closed for $k=1$.

In \cite{BanksElderWillis} the following properties of $G^{\propP{k}}$ are established for $G \leq \aut(T)$ where $T$ is locally finite; see \cite{ReidSmith} for trivial adjustments to the arguments so they continue to work in the non-locally finite case. Again let $T$ be a   (not necessarily locally finite)   tree, let $G \leq \aut(T)$, and fix $k \in \N$.
\begin{enumerate}
\item 
    $G^{\propP{k}}$ is a closed subgroup of $\Aut(T)$.
\item
    $G^{\propP{r}} = (G^{\propP{k}})^{\propP{r}}$ whenever $r \le k$.  In particular, $G^{\propP{k}}$ is $\propP{k}$-closed.
\item
    If $G$ is closed then $G = G^{\propP{1}}$ if and only if $G$ satisfies Tits' Independence Property $\propP{}$.  Consequently, to decide if $G$ satisfies Tits' Independence Property $\propP{}$ it is sufficient to check only that the definition of Property $\propP{}$ holds for paths in $T$ of length one (i.e. edges).
\end{enumerate}

The $\propP{k}$-closures of $G$ as $k \rightarrow \infty$ give a series of approximations to the action of $G$ on $T$. Moreover, the $\propP{k}$-closures converge to the closure of the action of $G$ in $\aut(T)$. Intuitively then the notion of $\propP{k}$-closure gives a tool for understanding all actions of groups on trees as `limits' of `universal groups' of `local actions' on $k$-balls. In our theory of local action diagrams we make precise this idea for $k=1$, creating an overarching theory for   $\propP{k}$-closed   groups when $k=1$. \\

Understanding   $\propP{k}$-closures   for $k>1$ presents significant challenges, but it is the obvious direction in which to generalise our theory of local action diagrams (see Section~\ref{Section:FutureWork}). Stephan Tornier in \cite{Tornier} has generalised the Burger--Mozes construction (for locally finite trees $T_d$) in a way that allows the local action on balls of a given radius $k \geq 1$ to be prescribed, where the Burger--Mozes groups $\U(F)$ then arise precisely when $k=1$. Even in this setting, Tornier's work highlights how the global behaviour of the resulting `universal' group is difficult to control---interactions between the prescribed local action on overlapping $k$-balls is the source of this complexity.

Let $\lcol$ be a legal colouring in the sense of Burger--Mozes (see Section~\ref{Section:BurgerMozes}) of the locally finite $d$-regular tree $T_d$.
Fix a subtree of the labelled tree $T_{d}$ arising as a ball of radius $k$ around some vertex, and denote this by $B_{d,k}$.
For each vertex $v \in VT_d$ recall that $B_v(k)$ is the $k$-ball around $v$, which is isomorphic (as a labelled tree) to $B_{d,k}$; denote this label-respecting isomorphism (which is unique) by $l_{v}^{k}:B_v(k)\to B_{d,k}$.
The {\it $k$-local action} of $g \in \aut(T_d)$ is then the graph isomorphism
\[\sigma_{k, v}(g): B_{d,k} \rightarrow B_{d,k} \quad \quad \sigma_{k, v}(g) = l_{gv}^{k}\circ g\circ (l_{v}^{k})^{-1}.\]
For any group $F \leq \aut(B_{d,k})$ we define $\U_k(F)$ as follows,
\[
 \U_{k}(F) := \{g \in \aut(T_{d}) : \sigma_{k, v}(g) \in F \quad \forall v\in VT_d\}.
\]

In \cite{Tornier} it is shown that $\U_{k}(F)$ is a closed, vertex transitive and compactly generated subgroup of $\aut(T_d)$, however $\U_{k}(F)$ can fail to have $k$-local action $F$. Despite this, the group retains a universal property: if $H \leq \aut(T_d)$ is locally transitive and contains an involutive inversion (that is, an involution $g$ that maps some $e \in AT_d$ to $\overline{e}$), and $F$ is the $k$-local action of $H$, then $\U_k(F)$ contains an $\aut(T_d)$-conjugate of $H$.

Precisely when $\U_{k}(F)$ has local action $F$ is characterised by a condition on $F$ that Tornier calls Condition (C); this condition can be thought of as a `compatibility' condition for the local action with itself as it interacts with an involutive automorphism. Another condition on $F$ called Condition (D) is sufficient for $\U_k(F)$ to be discrete, and when $F$ satisfies Condition (C) then (D) precisely characterises whether or not $\U_k(F)$ is discrete.\\

Finally we mention Caprace and De Medts' paper \cite{CapraceDeMedts} which contains an abundance of interesting results. The paper largely concerns compactly generated locally compact groups which act on locally finite trees and satisfy Tits’ independence property $\propP{}$. In other words, compactly generated locally compact $\propP{}$-closed actions on locally finite trees. 
The central theme of the paper is to consider the extent to which selected global properties of these groups (e.g.~having every proper open subgroup compact) are determined by the local action of a stabiliser of a vertex on its neighbours. The theory of local action diagrams shows that {\it all} properties of $\propP{}$-closed actions on trees are determined by these local actions, because every $\propP{}$-closed action is the universal group of a local action diagram.

%
%%
%%%
%%%%
%%%%%
%%%%%%
%%%%%%%
%%%%%%%%
%%%%%%%%%
%%%%%%%%%%
%%%%%%%%%%
%%%%%%%%%%
%%%%%%%%%%
\section{The theory of local action diagrams}

In this section we give an intuitive introduction to the various ingredients of our theory, and we describe the relationship between them. Everything here can be found in \cite{ReidSmith}.

%
%%
%%%
%%%%
%%%%%
%%%%%%
%%%%%%%
%%%%%%%%
%%%%%%%%%
%%%%%%%%%%
%%%%%%%%%%
%%%%%%%%%%
%%%%%%%%%%
\subsection{Local action diagrams}

At the heart of our theory is a combinatorial object we call a {\it local action diagram}. It plays the same role in our theory that a graph of groups plays in Bass--Serre Theory (see Section~\ref{Sec:BassSerre}).

\begin{definition}
A {\it local action diagram} is comprised of three things:
\begin{enumerate}
\item
    A nonempty connected graph $\Gamma$.
\item
    A nonempty set $X_a$ of colours for each arc $a \in A\Gamma$, such that the colour sets of distinct arcs are disjoint. For each vertex $v \in V\Gamma$ let $X_v := \bigsqcup_{a \in o\inv(v)}X_a$.
\item
	A closed group $G(v) \leq \sym(X_v)$ for each vertex $v \in V\Gamma$, such that the arc colour sets $X_a$ are the orbits of $G(v)$ on $X_v$.
\end{enumerate}
We denote such a local action diagram as $\Delta = (\Gamma,(X_a),(G(v)))$.
We call each $X_a$ the {\it colour set} of $a$, and each group $G(v)$ the {\it local action} at $v$.

  Notice how easy it is to construct a local action diagram. One chooses any nonempty connected graph and some groups (that have the specified orbit structure) to play the role of local actions; there are no further compatibility conditions governing whether or not various local actions can be combined on a local action diagram.  
Nevertheless, we will see they provide a complete description of all $\propP{}$-closed groups of automorphisms of trees.

Two local action diagrams are isomorphic if the graphs are isomorphic and their local actions are permutationally isomorphic, and these two types of isomorphisms are compatible. More precisely, if we have two local action diagrams $\Delta = (\Gamma,(X_a),(G(v)))$ and $\Delta' = (\Gamma',(X'_a),(G'(v)))$, an isomorphism between them consists of two things:

\begin{enumerate}
\item
    an isomorphism $\theta: \Gamma \rightarrow \Gamma'$ of graphs; and
\item
    bijections $\theta_v: X_v \rightarrow X'_{\theta(v)}$ for each vertex $v \in V\Gamma$, that restrict to a bijection $X_a \rightarrow X_{\theta(a)}$ for each $a \in o\inv(v)$, and such that $\theta_v G(v) \theta\inv_v = G'(\theta(v))$.
\end{enumerate}
\end{definition}

\begin{example} \label{Ex:TwoTypesOfLAD}
In Figures~\ref{Fig:MundaneLAD} and \ref{Fig:FancyLAD} we see examples   of local action diagrams.   For now they are nothing more than combinatorial objects, but we will soon see that they give rise to universal groups acting on trees. The second local action diagram shows how one can embed interesting permutation groups or topological groups as local actions: the automorphism group of the $3$-regular tree $T_3$ is a component of the local action of $v$, where $\aut(T_3)$ is acting on the colour set (i.e.~$VT_3$) of the arc $a$.

\begin{figure}
\begin{center}
\scalebox{0.85}{ % Usual scale command doesn't reduce fonts
\begin{tikzpicture}
\draw[rounded corners] (-5.7, -2) rectangle (5.7, 2) {};
%
% -= Vertex V
\node[Vertex] (V) at (-4,0) {};
    % Vertex label
    \node[TextNode] at ([shift=(180:0.3)]V) {$v$};
    % Local action
    \node[TextNode] at ([shift=(90:1.6)]V) {$G(v) = \langle (1 2 3)(8 9) \rangle$};
%
% -= Vertex W
\node[Vertex] (W) at (4,0) {};
    % Vertex label
	\node[TextNode] at ([shift=(0:0.3)]W) {$w$};
    % Local action
    \node[TextNode] at ([shift=(270:1.6)]W) {$G(w) = \langle (4 5)(6 7) \rangle$};
%
% -= Edge a [leaving v]
\draw [color=black, ->-, line width=0.8] (V) to [bend left=30] node[above] (colour_a) {$a$} (W);
    % Colour set
    \node[TextNode,align=left] at ([shift=(0:2)]colour_a) {$X_a = \{1,2,3\}$};
\draw [color=black, ->-, line width=0.8] (W) to [bend right=15]  node[below] (colour_aInv) {$\overline{a}$} (V);
    % Colour set
    \node[TextNode,align=left] at ([shift=(-8:1.4)]colour_aInv) {$X_{\overline{a}} = \{4,5\}$};
%
% -= Edge b [leaving v]
\draw [color=black, ->-, line width=0.8] (W) to [bend left=15]  node[above] (colour_bInv) {$\overline{b}$} (V);
    % Colour set
    \node[TextNode,align=left] at ([shift=(170:1.5)]colour_bInv) {$X_{\overline{b}} = \{6,7\}$};
\draw [color=black, ->-, line width=0.8] (V) to [bend right=30] node[below] (colour_b) {$b$} (W);
    % Colour set
    \node[TextNode,align=left] at ([shift=(180:2)]colour_b) {$X_{b} = \{8,9\}$};
\end{tikzpicture}
}
\end{center}
\caption{An example of a local action diagram.}
\label{Fig:MundaneLAD}
\end{figure}
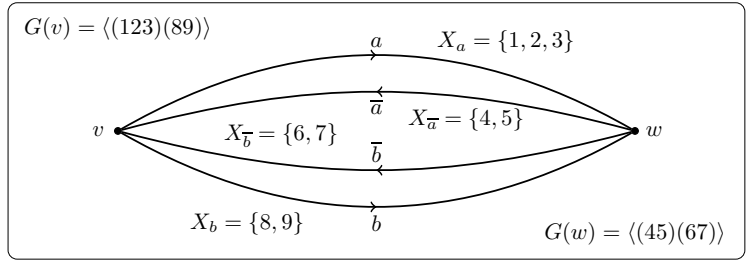

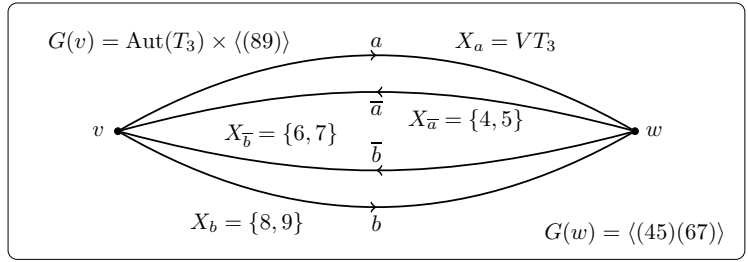
\begin{figure}
\begin{center}
\scalebox{0.85}{ % Usual scale command doesn't reduce fonts
\begin{tikzpicture}
\draw[rounded corners] (-5.7, -2) rectangle (5.7, 2) {};
%
% -= Vertex V
\node[Vertex] (V) at (-4,0) {};
    % Vertex label
    \node[TextNode] at ([shift=(180:0.3)]V) {$v$};
    % Local action
    \node[TextNode] at ([shift=(60:1.6)]V) {$G(v) = \aut(T_3) \times \langle (8 9) \rangle$};
%
% -= Vertex W
\node[Vertex] (W) at (4,0) {};
    % Vertex label
	\node[TextNode] at ([shift=(0:0.3)]W) {$w$};
    % Local action
    \node[TextNode] at ([shift=(270:1.6)]W) {$G(w) = \langle (4 5)(6 7) \rangle$};
%
% -= Edge a [leaving v]
\draw [color=black, ->-, line width=0.8] (V) to [bend left=30] node[above] (colour_a) {$a$} (W);
    % Colour set
    \node[TextNode,align=left] at ([shift=(0:2)]colour_a) {$X_a = VT_3$};
\draw [color=black, ->-, line width=0.8] (W) to [bend right=15]  node[below] (colour_aInv) {$\overline{a}$} (V);
    % Colour set
    \node[TextNode,align=left] at ([shift=(-8:1.4)]colour_aInv) {$X_{\overline{a}} = \{4,5\}$};
%
% -= Edge b [leaving v]
\draw [color=black, ->-, line width=0.8] (W) to [bend left=15]  node[above] (colour_bInv) {$\overline{b}$} (V);
    % Colour set
    \node[TextNode,align=left] at ([shift=(170:1.5)]colour_bInv) {$X_{\overline{b}} = \{6,7\}$};
\draw [color=black, ->-, line width=0.8] (V) to [bend right=30] node[below] (colour_b) {$b$} (W);
    % Colour set
    \node[TextNode,align=left] at ([shift=(180:2)]colour_b) {$X_{b} = \{8,9\}$};
\end{tikzpicture}
}
\end{center}
\caption{Interesting permutation groups can be embedded as local actions inside local action diagrams.}
\label{Fig:FancyLAD}
\end{figure}

\end{example}

%
%%
%%%
%%%%
%%%%%
%%%%%%
%%%%%%%
%%%%%%%%
%%%%%%%%%
%%%%%%%%%%
%%%%%%%%%%
%%%%%%%%%%
%%%%%%%%%%
\subsection{The associated local action diagram}
\label{Sec:AssocLad}

As we described in Section~\ref{sec:Assoc_graph}, for every group of automorphisms of a tree we can associate a graph of groups. The analogous situation arises here: for every group of automorphisms of a tree we can associate a local action diagram.

\begin{definition} \label{def:AssocLad}
Suppose $G$ is a group of automorphisms of a tree $T$. The {\it associated local action diagram} is a local action diagram $(\Gamma, (X_a), (G(v)))$ with the following parameters.
\begin{itemize}
\item
	For our connected graph we take $\Gamma$ to be the quotient graph $G \backslash T$, with $\pi$ denoting the natural quotient map.
\item
    For our arc colours we proceed as follows. The vertices of $\Gamma$ are orbits of $G$ on $VT$, so for each vertex $v \in V\Gamma$ we choose a representative vertex $v^*$ in $VT$ such that $\pi(v^*) = v$, and write $V^*$ for the set of all such representatives. Now the stabiliser $G_{v^*}$ permutes the arcs in $o\inv(v^*)$, and so we take $X_v := o\inv(v^*)$. The arcs in $\Gamma$ are the orbits of $G$ on $AT$, so this set $X_v$ breaks down into $G_{v^*}$-orbits as $X_a := \{b \in o\inv(v^*) : \pi(b) = a\}$ for each $a \in A\Gamma$ satisfying $o(a) = v$, and these sets $X_a$ become our arc colours.
\item
    For the local actions, take $G(v)$ to be the closure of the permutation group induced on $X_v$ by the stabiliser $G_{v^*}$, and note that the orbits of $G(v)$ and $G_{v^*}$ on $X_v$ coincide.
\end{itemize}
There are many choices for the associated local action diagram, but they are all isomorphic as local action diagrams.
\end{definition}

\begin{example} \label{Ex:AssocLAD}
Consider the diagram in Figure~\ref{Fig:AssocLAD}, which at the top shows part of the infinite $(2,4)$-biregular tree $T:= T_{2,4}$ (ignoring for now the decorations). Let $G$ be the closed subgroup of $\aut(T)$ with two orbits on vertices (indicated using black vertices and white vertices), and two orbits on edges (solid edges and dashed edges).

Thus the local action diagram has two vertices, say $v,w$, and two edges (each consisting of an arc in each direction), say $\{a, \overline{a}\}$ and $\{b, \overline{b}\}$. 
We pick vertex orbit representatives $v^*$ and $w^*$ in $VT$. The $T$-arcs leaving these vertices will become the arc colours in our local action diagram, so we name them as indicated: the two arcs leaving $v^*$ are named $1,2$ and the four arcs leaving $w^*$ are named $3$ to $6$.
Now $o\inv(v^*) = \{1,2\}$ and $G_{v^*}$ induces the trivial subgroup of $S_2$ on this set. Thus, we take $G(v) = \langle (1)(2) \rangle$. Meanwhile, $o\inv(w^*) = \{3,4,5,6\}$ and $G_{w^*}$ induces $\langle (3 4), (5 6) \rangle \leq \sym(\{3,4,5,6\})$ on $o\inv(w^*)$, so we take $G(w) = \langle (3 4), (5 6) \rangle$.

From these observations we obtain the local action diagram for $(T,G)$, shown at the bottom of Figure~\ref{Fig:AssocLAD}.

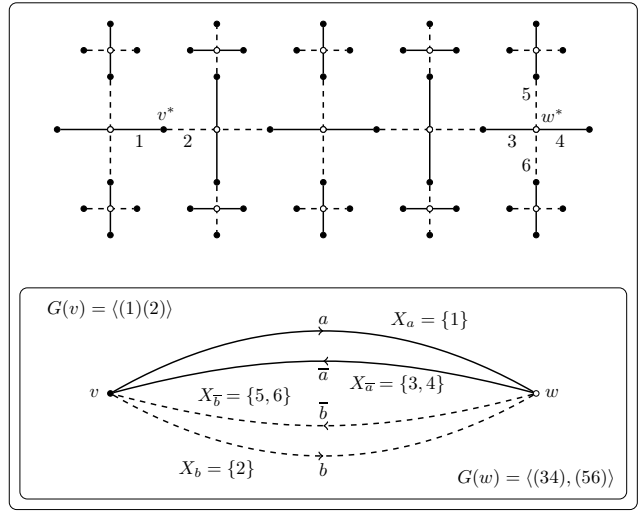
\begin{figure}
\begin{center}
\scalebox{0.7}{ % Usual scale command doesn't reduce fonts
\begin{tikzpicture}
%
% -= Tree
\begin{scope}
    \foreach \x in {-5,-3,...,5} {
        \node[Vertex] (VCen\x) at (\x,5) {};
    };
    \foreach \x in {-4,-2,...,4} {
        \node[OpenVertex] (VCen\x) at (\x,5) {};
        \foreach \y in {6} {
            \node[Vertex] (VMA\x) at (\x,\y) {};
            \node[OpenVertex] (VMU\x) at (\x,\y+0.5) {};
            \node[Vertex] (VTop\x) at (\x,\y+0.5+0.5) {};
            \node[Vertex] (VMUL\x) at (\x-0.5,\y+0.5) {};
            \node[Vertex] (VMUR\x) at (\x+0.5,\y+0.5) {};
        };
        \foreach \y in {4} {
            \node[Vertex] (VMB\x) at (\x,\y) {};
            \node[OpenVertex] (VML\x) at (\x,\y-0.5) {};
            \node[Vertex] (VBot\x) at (\x,\y-0.5-0.5) {};
            \node[Vertex] (VMLL\x) at (\x-0.5,\y-0.5) {};
            \node[Vertex] (VMLR\x) at (\x+0.5,\y-0.5) {};
        };
    };
    \foreach \x in {-4,0,4} {
        \draw [color=black, line width=0.8] (VMA\x) to (VTop\x);
        \draw [color=black, line width=0.8, dashed] (VMUL\x) to (VMUR\x);
        % Redraw white nodes so not covered by edge
        \node[OpenVertex] at (VMU\x) {};

        \draw [color=black, line width=0.8] (VMB\x) to (VBot\x);
        \draw [color=black, line width=0.8, dashed] (VMLL\x) to (VMLR\x);
        % Redraw white nodes so not covered by edge
        \node[OpenVertex] at (VML\x) {};

        \draw [color=black, line width=0.8, dashed] (VMA\x) to (VMB\x);
        % Redraw white nodes so not covered by edge
        \node[OpenVertex] at (VCen\x) {};
    };
    \foreach \x in {-2,2} {
        \draw [color=black, line width=0.8, dashed] (VMA\x) to (VTop\x);
        \draw [color=black, line width=0.8] (VMUL\x) to (VMUR\x);
        % Redraw white nodes so not covered by edge
        \node[OpenVertex] at (VMU\x) {};

        \draw [color=black, line width=0.8, dashed] (VMB\x) to (VBot\x);
        \draw [color=black, line width=0.8] (VMLL\x) to (VMLR\x);
        % Redraw white nodes so not covered by edge
        \node[OpenVertex] at (VML\x) {};

        \draw [color=black, line width=0.8] (VMA\x) to (VMB\x);
        % Redraw white nodes so not covered by edge
        \node[OpenVertex] at (VCen\x) {};
    };
    
    % Draw edges along central axis
    \draw [color=black, line width=0.8] (VCen-5)--(VCen-3);
    \draw [color=black, line width=0.8, dashed] (VCen-3)--(VCen-1);
    \draw [color=black, line width=0.8] (VCen-1)--(VCen1);
    \draw [color=black, line width=0.8, dashed] (VCen1)--(VCen3);
    \draw [color=black, line width=0.8] (VCen3)--(VCen5);
    % Redraw white vertices along centre
    \foreach \x in {-4,-2,...,4} {
        \node[OpenVertex] (VCen\x) at (\x,5) {};
    };
    
    % Label v* and w*
    \node[TextNode] at ([shift=(80:0.3)]VCen-3) {$v^*$};
    \node[TextNode] at ([shift=(40:0.4)]VCen4) {$w^*$};
    
    % Label arcs coming from v*
    \node[TextNode] at ([shift=(205:0.5)]VCen-3) {1};
    \node[TextNode] at ([shift=(-25:0.5)]VCen-3) {2};    
    % Label arcs coming from w*
    \node[TextNode] at ([shift=(205:0.5)]VCen4) {3};
    \node[TextNode] at ([shift=(-25:0.5)]VCen4) {4};
    \node[TextNode] at ([shift=(105:0.7)]VCen4) {5};
    \node[TextNode] at ([shift=(255:0.7)]VCen4) {6};
\end{scope}
%
% -= Local action diagram
\draw[rounded corners] (-5.7, -2) rectangle (5.7, 2) {};
%
% -= Vertex V
\node[Vertex] (V) at (-4,0) {};
    % Vertex label
    \node[TextNode] at ([shift=(180:0.3)]V) {$v$};
    % Local action
    \node[TextNode] at ([shift=(90:1.6)]V) {$G(v) = \langle (1)(2) \rangle$};
%
% -= Vertex W
\node[OpenVertex] (W) at (4,0) {};
    % Vertex label
	\node[TextNode] at ([shift=(0:0.3)]W) {$w$};
    % Local action
    \node[TextNode] at ([shift=(270:1.6)]W) {$G(w) = \langle (3 4), (5 6) \rangle$};
%
% -= Edge a [leaving v]
\draw [color=black, ->-, line width=0.8] (V) to [bend left=30] node[above] (colour_a) {$a$} (W);
    % Colour set
    \node[TextNode,align=left] at ([shift=(0:2)]colour_a) {$X_a = \{1\}$};
\draw [color=black, ->-, line width=0.8] (W) to [bend right=15]  node[below] (colour_aInv) {$\overline{a}$} (V);
    % Colour set
    \node[TextNode,align=left] at ([shift=(-8:1.4)]colour_aInv) {$X_{\overline{a}} = \{3,4\}$};
%
% -= Edge b [leaving v]
\draw [color=black, ->-, line width=0.8, dashed] (W) to [bend left=15]  node[above] (colour_bInv) {$\overline{b}$} (V);
    % Colour set
    \node[TextNode,align=left] at ([shift=(170:1.5)]colour_bInv) {$X_{\overline{b}} = \{5,6\}$};
\draw [color=black, ->-, line width=0.8, dashed] (V) to [bend right=30] node[below] (colour_b) {$b$} (W);
    % Colour set
    \node[TextNode,align=left] at ([shift=(180:2)]colour_b) {$X_{b} = \{2\}$};
    
% -= Final box
\draw[rounded corners] (-5.9, -2.2) rectangle (5.9, 7.4) {};
\end{tikzpicture}
}
\caption{Part of $T:=T_{2,4}$ with the arc and vertex orbits of $G \leq \aut(T)$ indicated (top) and the associated local action diagram for $(T,G)$ (bottom) described in Example~\ref{Ex:AssocLAD}.}
\label{Fig:AssocLAD}
\end{center}
\end{figure}
\end{example}

\begin{example} \label{Ex:AssocLADTwoLoops}
Let $d \geq 3$ be a finite or infinite cardinal and let $F \leq S_d$. Consider the Burger--Mozes group $\U(F) \leq \aut(T_d)$. Fix an arc $a \in AT_d$ and notice that the condition $\lcol(a) = \lcol(\overline{a})$ ensures that there is an element $g \in \aut(T_d)$ inverting $a$ (that is, mapping $a$ to $\overline{a}$ and vice versa) such that $\lcol(b) = \lcol(gb)$ for all $b \in AT_d$. For this element $g$ we have that $\sigma_{\lcol, v}(g)$ is trivial and therefore lies in $F$, for all $v \in VT_d$. Thus, $g \in \U(F)$.
We have seen then that for any arc $a \in AT_d$ we have that $a$ and $\overline{a}$ lie in the same orbit of $\U(F)$. Moreover, in Section~\ref{Section:BurgerMozes} we saw that $\U(F)$ is vertex transitive. Hence the associated local action diagram for $\U(F)$ consists of a single vertex with some loops, each of which is its own reverse. Each loop corresponds to an orbit of $F$.

For example, taking $F$ to be the group $\aut(T_3) \times S_3$ acting via the product action on $VT_3 \times \{1,2,3\}$, we have that $F$ has two orbits and therefore $\U(\aut(T_3) \times S_3)$ has the associated local action diagram shown in Figure~\ref{Fig:AssocLADTwoLoops}.

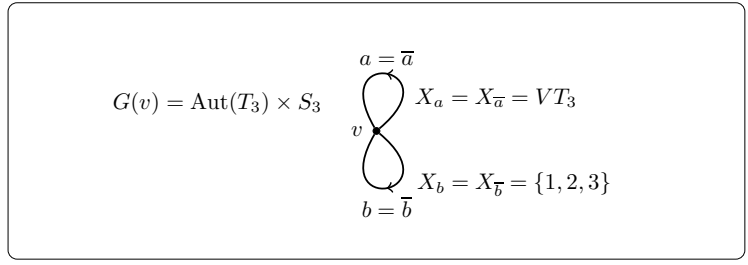
\begin{figure}
\begin{center}
\scalebox{0.85}{ % Usual scale command doesn't reduce fonts
\begin{tikzpicture}
\draw[rounded corners] (-5.7, -2) rectangle (5.7, 2) {};
%
% -= Vertex V
\node[Vertex] (V) at (0,0) {};
    % Vertex label
    \node[TextNode] at ([shift=(180:0.3)]V) {$v$};
    % Local action
    \node[TextNode] at ([shift=(170:2.5)]V) {$G(v) = \aut(T_3) \times S_3$};
% -= Loop a
    \draw[color=black, ->-, line width=0.8] (V) to  [out=40,in=120,distance=15mm] node[above] (colour_a) {$a = \overline{a}$} (V);
    % Colour set
    \node[TextNode,align=left] at ([shift=(-20:1.8)]colour_a) {$X_a = X_{\overline{a}} = VT_3$};
% -= Loop b
    \draw[color=black, ->-, line width=0.8] (V) to [out=-40,in=-120,distance=15mm] node[below] (colour_b) {$b = \overline{b}$} (V);
    % Colour set
    \node[TextNode,align=right] at ([shift=(10:2)]colour_b) {$X_b = X_{\overline{b}} = \{1,2,3\}$};    
\end{tikzpicture}
}
\caption{The associated local action diagram for Example~\ref{Ex:AssocLADTwoLoops}.}
\label{Fig:AssocLADTwoLoops}
\end{center}
\end{figure}
\end{example}

\begin{example} \label{Ex:SimpleLAD}
Let $d_1, d_2 > 1$ be finite or infinite cardinals, let $X_1, X_2$ be sets of cardinality $d_1, d_2$ respectively, and let $F_1 \leq \sym(X_1)$ and $F_2 \leq \sym(X_2)$. Consider the box product group $\U(F_1, F_2) \leq \aut(T)$, where $T$ is the biregular tree $T_{d_1, d_2}$.

In this general setting the local action diagram of $\U(F_1, F_2)$ can be significantly more complicated than those of the Burger--Mozes groups; despite this it is still tractable. Indeed, in \cite[Lemma 22]{SmithDuke} it is shown that the quotient graph $\U(F_1, F_2) \backslash T$ is the complete bipartite graph $K_{n_1,n_2}$, where $n_i$ is the number of orbits of $F_i$ on $X_i$ for $i=1,2$. In particular, between any two vertices in $K_{n_1,n_2}$ there is a single edge consisting of two arcs, one in each direction.

Identifying $\U(F_1, F_2) \backslash T$ and $K_{n_1,n_2}$ we can now form the associated local action diagram.
  Let $V_1, V_2$ be the two parts corresponding to the bipartition of $K_{n_1,n_2}$, with each part $V_i$ consisting of vertices with valency $n_i$ for $i=1,2$.  
Vertices in $V_i$ have local action $F_i$.   Each   arc $a$ from $V_1$ to $V_2$ represents an orbit $\Omega$ of $F_1$ and so we set $X_a := \Omega$. Similarly, each arc $b$ from $V_2$ to $V_1$ represents and orbit $\Omega'$ of $F_2$, and so we set $X_b := \Omega'$. We have a local action diagram.

In the special case where $F_1$ and $F_2$ are transitive, the local action diagram associated to $\U(F_1, F_2)$ is then two vertices connected by a single edge which consists of a pair of arcs, one in each direction. For example, the $2^{\aleph_0}$ pairwise nonisomorphic simple nondiscrete compactly generated \tdlc groups that were constructed in the proof of \cite[Theorem 38]{SmithDuke} were of the form $\U(\mathcal{O}_i, S_3)$, where each group $\mathcal{O}_i$ is viewed as a permutation group $\mathcal{O}_i \leq \sym(\mathcal{O}_i / H_i)$ as in Section~\ref{Section:BoxProduct}. The local action diagrams of these groups are shown in Figure~\ref{Fig:SimpleLAD}.

\begin{figure}
\begin{center}
\scalebox{0.85}{ % Usual scale command doesn't reduce fonts
\begin{tikzpicture}
\draw[rounded corners] (-5.7, -2) rectangle (5.7, 2) {};
%
% -= Vertex V
\node[Vertex] (V) at (-4,0) {};
    % Vertex label
    \node[TextNode] at ([shift=(180:0.3)]V) {$v$};
    % Local action
    \node[TextNode] at ([shift=(60:1.6)]V) {$G(v) = \mathcal{O}_i$};
%
% -= Vertex W
\node[Vertex] (W) at (4,0) {};
    % Vertex label
	\node[TextNode] at ([shift=(0:0.3)]W) {$w$};
    % Local action
    \node[TextNode] at ([shift=(270:1.6)]W) {$G(w) = S_3$};
%
% -= Edge a [leaving v]
\draw [color=black, ->-, line width=0.8] (V) to [bend left=15] node[above] (colour_a) {$a$} (W);
    % Colour set
    \node[TextNode,align=left] at ([shift=(0:2)]colour_a) {$X_a = \mathcal{O}_i / H_i$};
\draw [color=black, ->-, line width=0.8] (W) to [bend left=15]  node[below] (colour_aInv) {$\overline{a}$} (V);
    % Colour set
    \node[TextNode,align=left] at ([shift=(-10:1.5)]colour_aInv) {$X_{\overline{a}} = \{1,2,3\}$};
\end{tikzpicture}
}
\caption{Local action diagrams for the groups described in Example~\ref{Ex:SimpleLAD}.}
\label{Fig:SimpleLAD}
\end{center}
\end{figure}
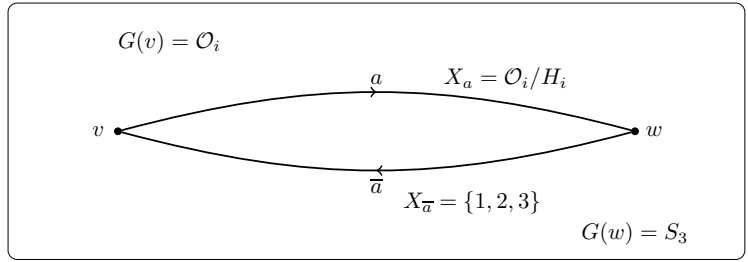
\end{example}

%
%%
%%%
%%%%
%%%%%
%%%%%%
%%%%%%%
%%%%%%%%
%%%%%%%%%
%%%%%%%%%%
%%%%%%%%%%
%%%%%%%%%%
%%%%%%%%%%
\subsection{Building a $\Delta$-tree} \label{Section:DeltaTree}

In Section~\ref{Sec:UnivCoveringBassSerre} we described how every graph of groups in Bass--Serre Theory gives rise to a universal cover that is a tree. An analogous situation arises here: every local action diagram gives rise to an arc-coloured tree we call a {\it $\Delta$-tree}. Intuitively this arc-coloured tree is built from taking `coloured walks' around the local action diagram. Shortly, we will construct a natural universal group from the automorphism group of this arc-coloured tree.

\begin{definition} \label{def:Delta_tree}
Let $\Delta = (\Gamma, (X_a), (G(v))$ be a local action diagram. A {\it $\Delta$-tree} is comprised of three things:
\begin{enumerate}
    \item 
        a tree $T$;
    \item
        a surjective graph homomorphism $\pi: T \rightarrow \Gamma$; and
    \item
        a colouring map $\lcol: AT \rightarrow \bigsqcup_{a \in A\Gamma}X_a$, such that for every vertex $v \in VT$, and every arc $a$ in $o\inv(\pi(v))$, the map $\lcol$ restricts to a bijection $\lcol_{v,a}$ from $\{b \in o\inv(v) : \pi(b) = a\}$ to $X_a$.
\end{enumerate}
We denote such a $\Delta$-tree as $\mathbf{T} = (T, \lcol, \pi)$. 
\end{definition}

In \cite[Lemma 3.5]{ReidSmith} it is shown that for any local action diagram $\Delta$ there exists a $\Delta$-tree, and moreover, any two $\Delta$-trees $(T,\pi,\lcol)$ and $(T',\pi',\lcol')$ are isomorphic in the following sense: there is a graph isomorphism $\alpha: T \rightarrow T'$ such that $\pi' \circ \alpha = \pi$.

The construction of a $\Delta$-tree is intuitive, but it requires careful notation to describe formally. We choose a root vertex $v_0 \in V\Gamma$. Then, for $v \in V\Gamma$ and $c \in X_v$, we define the {\it type} $p(c)$ of the colour $c$ to be the arc $a \in A\Gamma$ for which $c \in X_a$. A finite sequence of colours $(c_1,c_2,\dots,c_n)$ such that $o(p(c_{i+1})) = t(p(c_i))$ for all $1 \le i < n$ is called a {\it coloured path}, and we think of our vertices in $VT$ as being labelled by these coloured paths from the origin $v_0$.

Intuitively, we build $T$ inductively by taking all coloured paths around $\Gamma$ starting from $v_0$, but at each `step' we will have a choice to make about the reverse colour for our step. We can choose arbitrarily, subject to two constraints: the colour must be of the correct type, and once we have chosen a reverse colour we cannot pick it for the coloured arc corresponding to our next step. To make this precise, we need more notation.
If a vertex $v$ has label $(c_1,c_2,\dots,c_n)$ then we write the length of $v$ as $\ell(v) = n$, and we say any vertex with label $(c_1, c_2, \ldots, c_m)$ for $m \leq n$ is a {\it prefix} of $v$. If $v$ is a prefix of $w$ and $\ell(v) = \ell(w)-1$, then we write $v \ll w$. The {\it reverse label} of $v$ is $\overline{v} = (d_1,d_2,\dots,d_n)$, where each $d_i$ is a colour such that (i) $p(d_i) = \overline{p(c_i)}$, and (ii) whenever $v$ is a prefix of $w$, the label $\overline{v}$ is the corresponding prefix of $\overline{w}$.

We build the vertex set $VT$ of $T$ inductively, starting at a base vertex $()$. If we have already defined a vertex $v = (c_1,c_2,\dots,c_n)$ with $\overline{v} = (d_1,d_2,\dots,d_n)$, then we create new vertices $v_{+c_{n+1}} = (c_1,\dots,c_n,c_{n+1})$ by letting $c_{n+1}$ range over all colours satisfying $o(p(c_{n+1})) = t(p(c_n))$ and $c_{n+1} \neq d_n$. For each of these new vertices, choose $d_{n+1}$ arbitrarily from $X_{\overline{p(c_{n+1})}}$ and let $\overline{v_{+c_{n+1}}} =  (d_1,d_2,\dots,d_{n+1})$.

Now we have constructed the vertices of $T$, creating the remaining tree structure requires little effort. We define $AT_+ := \{(v,w) : v \ll w\}$ and $AT_- := \{(w,v) : (v,w) \in AT_+\}$ and $AT := AT_- \sqcup AT_+$, with functions $o, t$ and edge reversal defined in the obvious way. It is clear that from this we obtain a tree. Furthermore, $T$ is naturally arc-coloured by the following colouring map $\lcol$: for all $(v,w) \in AT_+$, take $\lcol(v,w)$ to be the last entry of $w$, and take $\lcol(w,v)$ to be the last entry of $\overline{w}$.

To define $\pi: T \rightarrow \Gamma$ we define $\pi$ of the base vertex $()$ to be our fixed root vertex $v_0 \in V\Gamma$, and then set $\pi(v)$ to be the vertex $t(p(c_n)) \in V\Gamma$, where $v \in VT$ has label $(c_1,\dots,c_n)$. For arcs in $a \in AT$, we set $\pi(a) = p(\lcol(a))$. It is now routine to verify that $\pi$ is a surjective graph homomorphism, and that $\lcol$ restricts to a bijection from $\{b \in o\inv(v) : \pi(b) = a\}$ to $X_a$. Whence, we have our $\Delta$-tree $\mathbf{T} = (T, \lcol, \pi)$.

\begin{remark}\label{Rem:TurnTintoDeltaTree}
Given a tree action $(T,G)$ and its associated local action diagram $\Delta$, we can equip $T$ with a natural colouring map $\lcol$ and projection map $\pi$ so that $(T, \lcol, \pi)$ is a $\Delta$-tree. Indeed, the associated local action diagram already has an appropriate map $\pi$, so we need only specify $\lcol$.
Using the notation of Definition~\ref{def:AssocLad} we choose, for each vertex $v \in VT$, an element $g_v \in G$ such that $g_v v \in V^*$.  Then $g_v$ induces a bijection from $o\inv(v)$ to $X_v$, because $X_v$ was defined to be $o\inv(v^*)$ for $v^* := g_v v \in V^*$. We then set $\lcol(b) := g_v b$ for all $b \in o\inv(v)$. It transpires that this colouring map satisfies the requirements of Definition~\ref{def:Delta_tree}, ensuring that $(T, \lcol, \pi)$ is indeed a $\Delta$-tree.
\end{remark}

%
%%
%%%
%%%%
%%%%%
%%%%%%
%%%%%%%
%%%%%%%%
%%%%%%%%%
%%%%%%%%%%
%%%%%%%%%%
%%%%%%%%%%
%%%%%%%%%%
\subsection{The universal group of a local action diagram}

In Section~\ref{Sec:BassSerre} we saw that every graph of groups in Bass--Serre Theory gives rise to a fundamental group that acts naturally on its universal cover. An analogous situation arises here: every local action diagram $\Delta$ gives rise to a universal group that acts naturally on a $\Delta$-tree.

\begin{definition}
Let $\Delta = (\Gamma, (X_a), (G(v)))$ be a local action diagram with $\Delta$-tree $\mathbf{T} = (T, \lcol, \pi)$.
An {\it automorphism} of $\mathbf{T}$ is a graph automorphism $\phi$ of the tree $T$ such that $\pi \circ \phi = \pi$. The group of all automorphisms of $\mathbf{T}$ is $\aut_{\pi}(T)$.
For $g \in \Aut_{\pi}(T)$ and $v \in VT$ we define the $\lcol$-local action of $g$ at $v$ as in Equation~\ref{eq:LocalAction}:
$\sigma_{\lcol,v}(g) := \lcol|_{o\inv(gv)} g \lcol|^{-1}_{o\inv(v)}$.
Again we note that $\sigma_{\lcol,v}(g) \in \sym(X_{\pi(v)})$ for all $v \in VT$.

The {\it universal group} of $\Delta$ and $\mathbf{T}$ is the group consisting of all $\mathbf{T}$-automorphisms whose local action at any $v \in VT$ always lies in the corresponding local action $G(\pi(v))$ of the local action diagram. Formally,
\[\U_{\mathbf{T}}(\Delta) := \{g \in \aut_{\pi}(T) : \sigma_{\lcol,v}(g) \in G(\pi(v)) \quad \text{ for all $v \in VT$}\}.\]
It transpires (see \cite[Theorem 3.12]{ReidSmith}) that   for a fixed $\Delta$,   different $\Delta$-trees give rise to the same universal group; that is, if $\mathbf{T}, \mathbf{T'}$ are $\Delta$-trees with underlying trees $T, T'$ respectively, then there is a graph isomorphism $\phi : T \rightarrow T'$ such that $\phi \U_{\mathbf{T}}(\Delta) \phi^{-1} = \U_{\mathbf{T'}}(\Delta)$. For this reason we typically omit the subscripts and speak of {\it the} universal group $\U(\Delta)$ of a local action diagram $\Delta$.
\end{definition}

\subsection{The correspondence theorem}
\label{Sec:LadCorrespondenceThm}

Recall that the $\propP{}$-closure of an action $(T, G)$ of a group $G$ on a tree $T$ is the smallest closed subgroup of $\aut(T)$ with Tits Independence Property $\propP{}$ that contains the action $(T,G)$. We denote this by $G^{\propP{}}$, and if $G = G^{\propP{}}$ then we say that $G$ is $\propP{}$-closed.

In Section~\ref{Sec:FundThmBassSerre} we described the Fundamental Theorem of Bass--Serre Theory, a correspondence theorem linking a tree action $(T,G)$ with the action of the fundamental group of its graph of groups on its universal cover. There is an analogous correspondence theorem for local action diagrams, linking the $\propP{}$-closure of a tree action $(T,G)$ with the action of the universal group of its local action diagram $\Delta$ on its $\Delta$-tree.

Moreover, in the Fundamental Theorem of Bass--Serre Theory we have that the fundamental group of a graph of groups $\mathbf{\Gamma}$ in its action on the universal cover of $\mathbf{\Gamma}$ (which is a tree) has $\mathbf{\Gamma}$ as its associated graph of groups. An analogous statement holds for local action diagrams: the universal group of a local action diagram $\Delta$ in its action on a $\Delta$-tree has an associated local action diagram that is isomorphic to $\Delta$.

\begin{theorem}[{\cite[Theorems 3.9 \& 3.10]{ReidSmith}}] \label{Thm:CorrespThmLad}
Let $G$ be a group and $T$ be a tree.

{\normalfont ($\circledast$)} Suppose $G$ acts on $T$ with associated local action diagram $\Delta$, universal group $\U(\Delta)$ and $\Delta$-tree $\mathbf{T}$.
Then the actions $(T, G^{\propP{}})$ and $(\mathbf{T}, \U(\Delta))$ can be identified.

{\normalfont ($\circledcirc$)} Suppose $\Delta$ is a local action diagram. Let $\mathbf{T}$ be a $\Delta$-tree and $\U(\Delta)$ its universal group. Then $\U(\Delta)$ is a $\propP{}$-closed group of automorphisms of the $\Delta$-tree $\mathbf{T}$ and the associated local action diagram of $(\mathbf{T}, \U(\Delta))$ is isomorphic to $\Delta$.
\end{theorem}

The first statement ($\circledast$) follows from the following observations. We have seen in Remark~\ref{Rem:TurnTintoDeltaTree} that 
we can equip $T$ with a colouring map $\lcol$ so that $T$ becomes a $\Delta$-tree $\mathbf{T'} = (T, \lcol, \pi)$. Using these we can construct the universal group $\U_{\mathbf{T'}}(\Delta) \leq \aut(T)$. In \cite[Theorem 3.10]{ReidSmith} we show that in fact $\U_{\mathbf{T'}}(\Delta) = G^{\propP{}}$. As noted previously, the possibly different $\Delta$-trees $\mathbf{T'}$ and $\mathbf{T}$ give rise to permutationally isomorphic universal groups, and via this relationship we can identify the actions of $(T, G^{\propP{}})$ and $(\mathbf{T}, \U(\Delta))$.
The second statement ($\circledcirc$) is \cite[Theorem 3.9]{ReidSmith}.

From this correspondence theorem, we obtain the following universal property for $\U(\Delta)$, which holds because for any group $G \leq \aut(T)$ with local action diagram $\Delta$ we have $G \leq G^{\propP{}} = \U(\Delta)$. This property clarifies the nature of the universal properties for the Burger--Mozes and box product groups, which both needed local transitivity to hold.

\begin{corollary}
Suppose $\Delta$ is a local action diagram, and form a $\Delta$-tree $\mathbf{T} = (T, \lcol, \pi)$. Then the universal group $\U(\Delta) \leq \aut(T)$ contains a permutationally isomorphic copy of every group $G \leq \aut(T)$ whose associated local action diagram is $\Delta$.
\end{corollary}

%
%%
%%%
%%%%
%%%%%
%%%%%%
%%%%%%%
%%%%%%%%
%%%%%%%%%
%%%%%%%%%%
%%%%%%%%%%
%%%%%%%%%%
%%%%%%%%%%
\section{A classification of groups with Tits' Independence Property $\propP{}$}

\subsection{The classification}

As we have already mentioned, Tits' Independence Property $\propP{}$ plays an important role in the theory of infinite groups, because it gives a general tool for constructing infinite nonlinear simple groups and because of its obvious importance to the subject of groups acting on trees. Our theory of local action diagrams gives a complete and highly usable description of all closed actions that have Tits' Independence Property $\propP{}$. From this one immediately obtains a classification of all groups that have Tits' Independence Property $\propP{}$.

In Theorem~\ref{Thm:CorrespThmLad} and the subsequent commentary, if we have a tree action $(T, G)$ where $G$ is closed with Tits' Independence Property $\propP{}$, then $G = G^{\propP{}}$ and therefore $(T, G)$ is equal to the universal group $\U_{\mathbf{T'}}(\Delta) \leq \aut (T)$, where $\Delta$ is the local action diagram of $(T,G)$ and $\mathbf{T'}$ is the $\Delta$-tree obtained from $T$ as described in Remark~\ref{Rem:TurnTintoDeltaTree}. On the other hand, if we have the universal group $\U(\Delta)$ of a local action diagram then it is a group of automorphisms of the underlying tree $T$ of a $\Delta$-tree, and the action $(T, \U(\Delta))$ is closed with Tits' Independence Property $\propP{}$. Thus, we have the following description of closed groups with $\propP{}$.

\begin{quote}
{\it Closed groups of automorphisms of trees with Tits' independence property $(\mathrm{P})$ are precisely the universal groups of local action diagrams.}
\end{quote}

In our paper \cite{ReidSmith}, we state the correspondence between $\propP{}$-closed actions and local action diagrams as follows.

\begin{theorem}[{\cite[Theorem 3.3]{ReidSmith}}] \label{thm:CorrespondenceThm}
There is a natural one-to-one correspondence between isomorphism classes of $\propP{}$-closed actions on trees and isomorphism classes of local action diagrams.
\end{theorem}

This correspondence is easy for us to state explicitly. For a local action diagram $\Delta$, we have a corresponding pair $(T,\U(\Delta))$, where $T$ is the underlying tree of a $\Delta$-tree. As previously discussed, different $\Delta$-trees give rise to permutationally isomorphic universal groups. Hence
the pair $(T,\U(\Delta))$ is unique up to isomorphisms.
Moreover, by construction we see that two isomorphic local action diagrams $\Delta$ and $\Delta'$ will produce isomorphic actions $(T,\U(\Delta))$ and $(T',\U(\Delta'))$. 
We have shown that there is a well-defined map $\theta$ from isomorphism class of actions $(T,G)$, where $T$ is a tree and $G$ is a $\propP{}$-closed group of automorphisms of $T$, to isomorphism classes of local action diagrams.
Our correspondence theorem (Theorem~\ref{Thm:CorrespThmLad}) shows that $\theta$ is a bijection.  Thus   $\theta$   is our claimed natural one-to-one correspondence.

  Recall that a permutation group $G \leq \sym(\Omega)$ and its closure in $\sym(\Omega)$ have the same orbits on all $n$-tuples of $\Omega$, for all $n \in \N$. Our complete description of all $\propP{}$-closed actions on trees immediately yields a useful classification of all (not just closed) actions $(T,G)$, where $T$ is a tree and $G$ has Tits' Independence Property $\propP{}$, whereby such actions $(T,G)$ are classified according to the isomorphism type of their associated local action diagram. In such a classification, every such action $(T,G)$ lies in precisely one class and the classification gives a complete description of the closure of $G$. Thus, we now have a deep understanding of all groups with Tits' Independence Property $\propP{}$.

\subsection{Some consequences}
\label{Sec:SomeConsequences}

There are many consequences to our description of all $\propP{}$-closed actions on trees. In this section we describe two (see \cite{ReidSmith} for more). 

The first is that any $\propP{}$-closed action on a tree can be described completely by drawing a local action diagram. This means that {\it all} properties of the action (e.g.~whether it is simple, geometrically dense, etc) can be read directly from the local action diagram. We explore this further in Section~\ref{Sec:SimplictyScoposEtc}.

The second is that for natural numbers $d, n$ there are only finitely many conjugacy classes of $\propP{}$-closed actions $(T,G)$ such that $T$ is locally finite of bounded valency $d$ and $G$ has at most $n$ vertex orbits. Indeed, any such action arises as the universal group of a local action diagram $\Delta = (\Gamma, (X_a), (G(v)))$ where $\Gamma$ has at most $n$ vertices and all groups $G(v)$ are finite permutation groups of degree at most $d$. In a sense, our theory reduces the study of such groups to questions about finite graphs and finite permutation groups.
In particular, for sensible choices for $d, n$ it would be possible to enumerate and describe (with the help of a computer) all such groups, for example by constructing all (finitely many) possible local action diagrams and then determining which are isomorphic.

In \cite[\S 7.1]{ReidSmith} we consider a special case of this: $\propP{}$-closed actions $(T_d,G)$ where $T_d$ is the $d$-regular tree and $G$ is vertex-transitive. By our classification theorem we know that up to conjugacy, there are only finitely many such actions. To determine them, we define an {\it orbit pairing} for $H \leq S_d$ to be a permutation of the set $H\backslash [d]$ of $H$-orbits whose square is the identity,   where $[d]$ here denotes the set $\{1,\ldots,d\}$.  
We consider pairs $(H, r)$, where $H \leq S_d$ and $r$ is an orbit pairing for $H$, and say two pairs $(H_1, r_1)$ and $(H_2, r_2)$ are equivalent whenever there exists 
$g \in S_d$ such that $gH_1g\inv = H_2$ and the map $g': H_1 \backslash [d] \rightarrow H_2 \backslash [d]$ induced by $g$ satisfies $g'r_1 = r_2g'$. The (finitely many) $\Aut(T_d)$-conjugacy classes for $\propP{}$-closed vertex transitive actions $(T_d,G)$ are in one-to-one correspondence with the set of equivalence classes of pairs $(H,r)$.

Since each $G$ is vertex transitive, its associated local action diagram $\Delta = (\Gamma, (X_a), (G(v)))$ has only one vertex $v_0$, together with some loops that may or may not be their own reverse. For each pair $(H, r)$, the group $H$ is the group $G(v_0)$; the arcs of $\Delta$ correspond to orbits of $H$; and the edges of $\Delta$ correspond to orbits of $r$ on $H \backslash [d]$. Recall that the vertices, arcs and edges of $\Delta$ correspond to, respectively, the vertex-, arc- and edge-orbits of $\U(\Delta)$ on $T$. Those orbits of $H$ that are fixed by $r$ correspond to the arc-orbits of $\U(\Delta)$ on arcs that are reversed by some element in $\U(\Delta)$. The equivalence classes of these pairs $(H,r)$ give rise to all isomorphism classes of local action diagrams of $\propP{}$-closed vertex transitive actions $(T_d,G)$, and we can enumerate these pairs for reasonable values of $d$.

In \cite[\S 7]{ReidSmith} we use this method to classify all such actions for $0 \leq d \leq 5$. The appendix to \cite{ReidSmith} is written by Stephan Tornier and contains a GAP (\cite{GAP4}) implementation that can perform this classification for values of $d$ greater than $5$. Even for $d=3$ we find examples of such actions that do not arise as Burger--Mozes groups; for larger values of $d$ the GAP implementation shows that the conjugacy classes of vertex transitive $\propP{}$-closed actions on $T_d$ that do not arise as Burger--Mozes groups grows rapidly.

%
%%
%%%
%%%%
%%%%%
%%%%%%
%%%%%%%
%%%%%%%%
%%%%%%%%%
%%%%%%%%%%
%%%%%%%%%%
%%%%%%%%%%
%%%%%%%%%%
\section{Reading simplicity from a local action diagram}
\label{Sec:SimplictyScoposEtc}

Recall Tits' result, Theorem~\ref{theorem:Tits4_5}: If $G \leq \aut(T)$ has property $\propP{}$ then the subgroup $G^+$ generated by arc stabilisers is trivial or simple whenever the action of $G$ is geometrically dense (i.e. $G$ leaves invariant no nonempty proper subtree of $T$ and $G$ fixes no end of $T$). Closed groups $G \leq \aut(T)$ with property $\propP{}$ are completely described by their local action diagrams, so it is not surprising that the simplicity of $G^+$ can be read directly from the local action diagram. What is surprising, however, is how easily this can be done. It happens that the invariant subtrees and fixed ends of a tree action $(T,G)$ correspond to combinatorial features of the local action diagram that we call {\it strongly confluent partial orientations} (or \scpos). Again we note that $T$ does not need to be locally finite.

\begin{definition}
A {\it strongly confluent partial orientation} (or {\it \scpo}) in a local action diagram $\Delta = (\Gamma,(X_a),(G(v)))$ is a subset $O$ of $A\Gamma$ satisfying:
\begin{enumerate}
\item
    For all $a \in O$ we have $\overline{a} \not\in O$ and $|X_a|=1$;
\item \label{item:scopo}
    For all $a \in O$ we have that $t\inv(o(a)) \setminus \{\overline{a}\} \subseteq O$.
\end{enumerate}
The empty set is always a \scpo; if the empty \scpo is the only \scpo of $\Delta$ then $\Delta$ is said to be {\it irreducible}.
\end{definition}

In \cite[Theorem 1.4]{ReidSmith} it is noted that the invariant subtrees and fixed ends of $(T, G)$ correspond to \scpos of the local action diagram $\Delta$. Under this natural correspondence, the empty \scpo (which always exists) corresponds to $T$ (which is always invariant under $G$). Thus $(T,G)$ being geometrically dense is equivalent to $\Delta$ being irreducible.

We can completely characterise all types of \scpos that correspond to invariant subtrees and fixed ends of faithful  actions $(T,G)$ with Property $\propP{}$. There are four types of these \scpos, and we call them a {\it stray leaf}, a {\it focal cycle}, a {\it horocyclic end} and a {\it stray half-tree}. Before describing these types of \scpos, we note the following.

\begin{theorem}[{\cite[Corollary 1.5]{ReidSmith}}] If $G \leq \Aut(T)$ has Tits' Property $\propP{}$ and local action diagram $\Delta$, then the following are equivalent:
\begin{enumerate}
    \item
        $G$ is geometrically dense;
    \item \label{item:NoFocalCycle}
        $\Delta$ is not a focal cycle and has no stray half-trees, no horocyclic ends and no stray leaves.
\end{enumerate} 
In particular, if \ref{item:NoFocalCycle} holds then $G^+$ is abstractly simple or trivial.
\end{theorem}
In this theorem we have a complete characterisation in the local action diagram of when Tits' Theorem can be applied to an action with Tits' Property $\propP{}$.

\begin{definition} Let $\Delta = (\Gamma,(X_a),(G(v)))$ be a local action diagram.
The following are the aforementioned \scpos that correspond to fixed ends and proper invariant subtrees.
\begin{itemize}
    \item
        If $\Gamma$ is a finite cycle with a cyclic orientation $O$, such that for all $a \in O$ we have $|X_a|=1$, then we say that $\Delta$ is a {\it focal cycle}.
    \item
        A {\it stray leaf} of $\Delta$ is a leaf $v$ of $\Gamma$ such that $|X_v| = 1$ (or equivalently such that $G(v)$ is trivial).
    \item
        A {\it horocyclic end} of $\Delta$ occurs only when $\Gamma$ is a tree. It is an end $\xi$ of the tree $\Gamma$ such that all the arcs $a \in A\Gamma$ that are directed towards $\xi$ satisfy $|X_a|=1$.
    \item
        If $\Gamma \setminus \{a, \overline{a}\}$ is not connected and $\Gamma_a$ is the connected component containing $t(a)$, then $\Gamma_a$ is a {\it stray half-tree} of $\Delta$ whenever $\Gamma_a$ is a tree that contains no leaves of $\Gamma$, and moreover within $\Gamma_a$ all arcs $b$ orientated towards $t(a)$ satisfy $|X_{b}|=1$.
\end{itemize}
\end{definition}

This characterisation of \scpos allows us to quickly determine whether or not $\Delta$ is irreducible in the frequently encountered case where $\Gamma$ is a finite graph that is not a cycle: $\Delta$ is irreducible if and only if $\Delta$ has no stray leaves. Thus, in this situation, if $\Delta$ has no stray leaves then $G^+$ is simple or trivial.

Determining whether or not $G^+$ is trivial is easy to detect in the local action diagram.

\begin{proposition}[{\cite[Lemma 5.11]{ReidSmith}}] If $T$ is a tree and $G \leq \aut(T)$ with associated local action diagram $\Delta = (\Gamma,(X_v),(G(v)))$, then $G^+$ is trivial if and only if $G(v)$ acts freely (i.e.~semiregularly) on $X_v$ for all $v \in V\Gamma$.
\end{proposition}

Our discussion so far concerned the simplicity of $G^+$. However, in various naturally occurring situations (see for example property \ref{item:BoxSimple} of the box product construction in Section~\ref{Section:BoxProduct}) it is the case that $G$ itself is simple. Regarding this, we have an almost   complete   characterisation of simplicity for faithful $\propP{}$-closed actions. It is almost   complete   because we must exclude some degenerate cases and we must ensure that the action results in a closed induced action on any invariant subtree.

\begin{definition}
A group $G \leq \Aut(T)$ is {\it strongly closed} if for every $G$-invariant subtree $T'$ of $T$, the induced action of $G$ on $T'$ is closed.
\end{definition}

Being strongly closed is easily achieved. For example, in \cite[Corollary 6.4]{ReidSmith}, we show that a locally compact $\propP{}$-closed subgroup of $\Aut(T)$ that acts with translation is always strongly closed.

\begin{theorem}[{\cite[Theorem 1.8]{ReidSmith}}]
If $(T,G)$ is a faithful $\propP{}$-closed and strongly closed action on a tree $T$, then the following are equivalent:
\begin{enumerate}
\item
    $G$ is nondiscrete, abstractly simple, and acts with translation.
\item
    There exists an infinite $G$-invariant subtree $T'$ of $T$ (not necessarily proper) on which $G$ acts faithfully. Furthermore, if the associated local action diagram of $(T',G)$ is $\Delta = (\Gamma,(X_a),(G(v)))$, then 
            $\Gamma$ is a tree;
            $\Delta$ is irreducible;
            all the groups $G(v)$ are closed and generated by point stabilisers; and
            at least one of the groups $G(v)$ is nontrivial.
\end{enumerate}
\end{theorem}

%
%%
%%%
%%%%
%%%%%
%%%%%%
%%%%%%%
%%%%%%%%
%%%%%%%%%
%%%%%%%%%%
%%%%%%%%%%
%%%%%%%%%%
%%%%%%%%%%
\section{Topological properties of universal groups}

In this section we survey a selection of results in \cite{ReidSmith} concerning various topological properties of $\propP{}$-closed subgroups of $\Aut(T)$. All statements are with respect to the permutation topology. Of course the topological properties of $\propP{}$-closed subgroups of $\Aut(T)$ when $T$ is locally finite are well-understood. The novelty in our results is that (i) we make no assumptions about local finiteness, and (ii) our results are typically concerned with deducing `global' topological statements from `local' properties that can be found in the local action diagram.

To characterise local compactness and compact generation of `nondegenerate' $\propP{}$-closed groups via their local action diagrams, we first need to define a combinatorial feature of local action diagrams called a {\it cotree}.
Let $\Gamma$ be a connected graph. Directed paths $(v_0,\dots,v_n)$ in $\Gamma$ are called {\it nonbacktracking} if $n = 0, 1$ or for $n \geq 2$ we have that $v_i \not = v_{i+2}$ for all $0 \leq i \leq n-2$.
Given an induced subgraph $\Gamma'$ of $\Gamma$, a finite directed nonbacktracking path $(v_0,\dots,v_n)$ in $\Gamma$ such that $v_n \in V\Gamma'$ and $v_i \not\in V\Gamma'$ for $i < n$, is called a {\it projecting path} from $v_0 \in V\Gamma$ to $\Gamma'$.
We say that $\Gamma'$ is a {\it cotree of $\Gamma$} if it is nonempty and for all $v \in V\Gamma \setminus V\Gamma'$ there is precisely one projecting path $(v_0,\dots,v_n)$ from $v$ to $\Gamma'$ (in particular this means that multiple arcs in $\Gamma$ from any $v_i$ to $v_{i+1}$ are not permitted).
Note that if $\Gamma$ is a connected graph that is not a tree then there is always a (unique) smallest cotree of $\Gamma$
and cotrees are connected induced subgraphs that contain this smallest cotree.
Given a cotree $\Gamma'$ of $\Gamma$, there is a \scpo $O_{\Gamma'}$ associated with $\Gamma'$, consisting of all arcs $a$ satisfying (i) $o(a) \not\in V\Gamma'$ and (ii) $a$ lies on the projecting path from $o(a)$ to $\Gamma$.

By excluding some degenerate cases, we can characterise local compactness of $\propP{}$-closed actions via the local action diagram as follows.

\begin{proposition}[{\cite[Proposition 6.3]{ReidSmith}}]
Suppose $T$ is a tree and $G \leq \Aut(T)$. Let $\Delta = (\Gamma,(X_a),(G(v)))$ be the local action diagram for $(T,G)$ and let $\Gamma'$ be the unique smallest cotree of $\Delta$. Suppose further that there is a unique minimal $G$-invariant subtree $T'$ of $T$ that has at least $3$ vertices. Then the following are equivalent.
\begin{enumerate}
\item $G^{\propP{}}$ is locally compact.
\item For all $a \in AT'$, the arc stabiliser $(G^{\propP{}})_a$ is compact.
\item For all $a \in A\Gamma$ such that $\overline{a} \not\in O_{\Gamma'}$, and for all $x \in X_a$, the orbits of the stabiliser $(G(o(a)))_x$ in its action on $X_v$ are finite.
\end{enumerate}
\end{proposition}

By excluding some degenerate cases, we can characterise the compact generation of $\propP{}$-closed actions via their local action diagrams as follows.

\begin{proposition}[{\cite[Proposition 6.5]{ReidSmith}}]
Let $T$ be a tree and suppose $G \leq \Aut(T)$ is closed with all vertex-orbits having unbounded diameter. Let $\Delta = (\Gamma,(X_a),(G(v)))$ be the local action diagram of $(T,G)$. If some arc stabiliser in $G$ is compact, then $G$ and $G^{\propP{}}$ are locally compact, and the following are equivalent.
\begin{enumerate}
\item $G$ is compactly generated;
\item $G^{\propP{}}$ is compactly generated;
\item there is a unique smallest $G$-invariant subtree $T'$ such that $G$ has finitely many orbits on $VT' \sqcup AT'$ and $G_v$ is compactly generated for each $v \in VT'$;
\item there is a unique smallest cotree $\Gamma'$ of $\Delta$ such that $\Gamma'$ is finite and $G(v)$ is compactly generated for each $v \in V\Gamma'$.
\end{enumerate}
\end{proposition}

Recall the class $\mathscr{S}$ of nondiscrete, topologically simple, compactly generated, locally compact groups. Let $\mathscr{S}_{td}$ be the class of totally disconnected groups in $\mathscr{S}$.
For constructing groups in the class $\mathscr{S}_{td}$, we are typically interested in the situation where the universal group is compactly generated, locally compact, and acts geometrically densely on its associated tree (allowing us to deduce simplicity via Tits' Theorem). For this situation we have the following.

\begin{theorem}[{\cite[Theorem 1.9]{ReidSmith}}]
Suppose $\Delta = (\Gamma,(X_a),(G(v)))$ is a local action diagram. Then $\U(\Delta)$ is compactly generated, locally compact, and acts geometrically densely on its associated tree if and only if $\Delta$ is irreducible, $\Gamma$ is finite, and all groups $G(v)$ are subdegree-finite and compactly generated.
\end{theorem}

The following result is particularly useful, since it allows us to construct local action diagrams that immediately yield groups in $\mathscr{S}_{td}$. Again for faithful $\propP{}$-closed actions it is an almost perfect characterisation of membership of the class $\mathscr{S}_{td}$.

\begin{theorem}[{\cite[Corollary 1.10]{ReidSmith}}]
If $(T,G)$ is a faithful $\propP{}$-closed and strongly closed action on a tree $T$, then the following are equivalent:
\begin{enumerate}
\item
    $G \in \mathscr{S}_{td}$ and $G$ fixes no vertex of $T$.
\item \label{item:lad_description:comp_gen}
    There exists a unique smallest $G$-invariant subtree $T'$ of $T$ (not necessarily proper) on which $G$ acts faithfully.  Furthermore, if $\Delta = (\Gamma,(X_a),(G(v)))$ is the associated local action diagram of $(T',G)$, then $\Gamma$ is a finite tree; all of the groups $G(v)$ are closed, compactly generated, subdegree-finite and generated by point stabilisers; and for every leaf $v$ of $\Gamma$ the group $G(v)$ is nontrivial.
\end{enumerate}
\end{theorem}

To conclude this section we give a theorem that establishes an entirely different universal property of the groups $\U(\Delta)$ within the class $\mathscr{S}_{td}$.

\begin{theorem}[{\cite[Theorem 1.14]{ReidSmith}}]
Let $G_1,\dots,G_n$ be a finite list of nontrivial compactly generated \tdlc groups, such that for each $G_i$ there is a compact open subgroup $U_i$ such that $G_i = \langle g U_i g\inv : g \in G_i \rangle$ and $\bigcap_{g \in G_i}gU_ig\inv$ is trivial. For example, we can take $G_i \in \mathscr{S}_{td}$ and $U_i$ to be any compact open subgroup.  Then there exists $\U(\Delta) \in \mathscr{S}_{td}$ acting continuously on a countable tree $T$, vertex stabilisers $O_1,\dots,O_n$ of $\U(\Delta)$ and compact normal subgroups $K_i$ of $O_i$, such that $O_i \cong K_i \rtimes G_i$ for $1 \le i \le n$.
\end{theorem}

%
%%
%%%
%%%%
%%%%%
%%%%%%
%%%%%%%
%%%%%%%%
%%%%%%%%%
%%%%%%%%%%
%%%%%%%%%%
%%%%%%%%%%
%%%%%%%%%%
\section{Project ideas and open problems}
\label{Section:FutureWork}

There are a number of directions in which our theory of local action diagrams can be generalised, and many areas where it can be applied; we describe some of them here.

\begin{project}
The application of local action diagrams to understand the automorphism groups of various types of infinitely ended graphs.
\end{project}
A connected, locally finite graph $\Gamma$ with infinitely many thin ends (i.e.~ends that do not contain infinitely many disjoint rays) and only countably many thick ends (i.e.~ends that contain infinitely many disjoint rays) is known by work of Carsten Thomassen and Wolfgang Woess (see \cite{ThomassenWoess}) to resemble a tree in some precise way. Thomassen and Woess' work relies heavily on Dicks and Dunwoody's theory of structure trees (see \cite{DicksDunwoody}), and it is via this theory that we can see the `tree-like' nature of $\Gamma$. See \cite{RoggiSurvey} for an `accessible' introduction to these ideas.

Now the automorphism group $\Aut(\Gamma)$ will act on the structure tree $T$ of $\Gamma$, and this action can be studied via local action diagrams. We give an example of this in \cite[\S7.2]{ReidSmith}, where we 
find a complete description of all automorphism groups of simple, nontrivial, vertex-transitive graphs with vertex connectivity one: they are precisely the universal groups of a certain type of local action diagram.

The automorphism groups of many other classes of graphs could be understood in this way.

\begin{project}
Generalise the ideas behind local action diagrams to better understand $\propP{k}$-closures of groups acting on trees.
\end{project}

This project idea appears in our paper (\cite[\S8 Question 2]{ReidSmith}). Such a project would be a significant undertaking, given the complexities of Tornier's generalisation of the Burger--Mozes groups (see Section~\ref{Sec:FurtherGens}). Nevertheless, it is our opinion that a usable theory could be developed, perhaps using a modified version of the local action action diagram built around $k$-arcs rather than arcs.

\begin{project}
Find further examples of locally determined global properties of $(T,G)$.
\label{Q:find_more_lad_features}
\end{project}

This is \cite[\S8 Question 5]{ReidSmith}. A long term research theme could be built around continuing to find global properties of $G \leq \Aut(T)$ that are perfectly characterised by the associated local action diagram. Such properties could be found by looking for global properties of $G$ that are characterised by properties of $G^{\propP{}}$ in its action on $T$. In \cite{ReidSmith} we call such properties {\it locally determined global properties of $(T,G)$}. As we have seen, geometrically dense actions are a locally determined global property. An important measure of success here will be that the locally determined global properties should be expressed in terms of features of the local action diagram.

\begin{project}
Write software to search local action diagrams for known features that characterise the locally determined global properties of actions $(T,G)$.
\end{project}

As we saw in \ref{Sec:SomeConsequences}, for given natural numbers $d,n$ there are only finitely many isomorphism classes of local action diagrams for actions $(T,G)$ where $G$ has at most $n$ vertex orbits on $T$, and $T$ is locally finite with every valency bounded by $d$. Constructing all possible local action diagrams for a given pair $(d,n)$ is then possible for reasonable choices of $(d,n)$, since these are finite graphs decorated with finite groups and finite colour sets. Algorithms can be created to search these constructions for known locally determined global properties. In this way, we could obtain lists of $\aut(T_d)$-conjugacy classes of, for example, vertex transitive actions on $T_d$ that have some locally determined global property.
As new properties come to light via Project~\ref{Q:find_more_lad_features}, more software solutions will be needed.

In \cite[\S8 Question 2]{ReidSmith} we give a somewhat related project, asking for the asymptotics of the number $N_d$ of conjugacy classes of vertex-transitive $\propP{}$-closed subgroups of $\Aut(T_d)$ as a function of $d$.

\begin{project}
Create analogies of the constructions in Section~\ref{Sec:FurtherGens}, for universal groups of local action diagrams.
\end{project}

In Section~\ref{Sec:FurtherGens} we outlined some generalisations of either the Burger--Mozes groups or the box product construction. For each of these generalisations there are natural ways to modify their definitions so that they can be applied to local action diagrams. For example, inspired by Le Boudec groups one could allow finitely many (or boundedly many) singularities for local actions in local action diagrams, and the local action at these singularities could be restricted to give an analogy to Le Boudec's restricted universal groups.

For all of these generalisations, it would be interesting to understand what permutational and topological properties the resulting `modified universal' groups possess.

\begin{furtherprojects} As this is an introductory note, we have omitted several aspects of the theory (e.g.~local subaction diagrams as a tool to better understand subgroups of universal groups of local action diagrams). There are several interesting questions arising from these omitted topics. We refer the interested reader to \cite[\S8]{ReidSmith}, where further open questions and project ideas are given.
\end{furtherprojects}


\begin{thebibliography}{8}

\bibitem{BCGM}
    Uri Bader, Pierre-Emmanuel Caprace, Tsachik Gelander and Shahar Mozes,
    Simple groups without lattices.
    {\it Bull.~Lond.~Math.~Soc.}
    44 (2012), 55--67.

\bibitem{BanksElderWillis}
    Christopher Banks, Murray Elder and George A.~Willis,
    Simple groups of automorphisms of trees determined by their actions on finite subtrees. 
    {\it J. Group Theory} 18 (2015), 235--261.

\bibitem{JensTom}
    Jens Bossaert and Tom De Medts,
    Topological and algebraic properties of universal groups for right-angled buildings.
    {\it Forum Math.} 33 (2021), 867--888.

\bibitem{BurgerMozes}
	Marc Burger and Shahar Mozes,
	Groups acting on trees: from local to global structure.
	{\it Publ. Math. IH\'ES}
	92 (2000), 113--150.

\bibitem {CapraceMonod} Pierre-Emmanuel Caprace and Nicolas Monod,
	Decomposing locally compact groups into simple pieces.
	{\it Math.~Proc.~Camb.~Phil.~Soc.}
	150 (2011), 97--128.
	
\bibitem {CapraceSimple} Pierre-Emmanuel Caprace,
	Non-discrete simple locally compact groups. In: 
	{\it European congress of mathematics. Proceedings of the 7th ECM (7ECM) congress, Berlin, Germany, July 18--22, 2016,}
	(European Mathematical Society (EMS), Z\"{u}rich, 2018), pp.~333--354.

\bibitem {CapraceDeMedts} Pierre-Emmanuel Caprace and Tom De Medts,
    Simple locally compact groups acting on trees and their germs of automorphisms.
    {\it Transformation Groups}, 16 (2011) 375--411.

\bibitem {Cam53}
	Ruth Camm,
	Simple free products.
	{\it J. London Math. Soc.}
	28 (1953), 66--76.

\bibitem{DicksDunwoody}
    Warren Dicks and M.~J.~Dunwoody,
    {\it Groups acting on graphs}.
    Cambridge Studies in Advanced Mathematics 
    (Cambridge University Press, 1989).

\bibitem{GAP4}
    The GAP~Group, GAP -- Groups, Algorithms, and Programming, Version 4.8.7 (2017), \url{http://www.gap-system.org}

\bibitem{goldschmidt}
	David M.~Goldschmidt,
	Automorphisms of trivalent graphs.
	{\it Ann.~Math.} 
	111 (1980), 377--406.

\bibitem{JungWatkins}
    H.~A.~Jung and M.~E.~Watkins,
    On the structure of infinite vertex-transitive graphs.
    {\it Discrete Math.} 
    18 (1977), 45--53.
	
\bibitem{LeBoudec}
    Adrien Le Boudec,
    Groups acting on trees with almost prescribed local action.
    {\it Commentarii Mathematici Helvetici} 91 (2016), 253--293.

\bibitem{Lederle}
    Waltraud Lederle,
	Coloured Neretin groups.
	{\it Groups Geom.~Dyn.} 13 (2019), 467--510.

\bibitem{Moller:Primitive}
    R\"{o}gnvaldur G. M\"{o}ller,
    Primitivity and ends of graphs.
    {\it Combinatorica}
    14 (1994), 477--484.

\bibitem{RoggiSurvey}
    R\"{o}gnvaldur G. M\"{o}ller,
    Groups acting on locally finite graphs—a survey of the infinitely ended case.
    In: {\it Groups `93 Galway/St. Andrews. Proceedings of the international conference, Galway, Ireland, August 1--14, 1993. Volume 2.}
    Lond.~Math.~Soc.~Lect.~Note Ser.~212,
    (Cambridge University Press, 1995), pp.~426--456.

\bibitem{Moller:PermTDLC}
    R\"{o}gnvaldur G. M\"{o}ller,
    Structure theory of totally disconnected locally compact groups via graphs and permutations.
    {\it Canad.~J.~Math.}
    54 (2002), 795--827.

\bibitem{potocnik_spiga_19}
	Primo\v{z} Poto\v{c}nik and Pablo Spiga,
	Lifting a prescribed group of automorphisms of graphs.
	{\it Proc.~Am.~Math.~Soc.}
	147 (2019), 3787--3796.

\bibitem{olshanski}   A.~Yu.~Ol'shanski{\u\i},  
	{\em Geometry of defining relations in groups}.
	Mathematics and its Applications (Soviet Series) 70
	(Kluwer Acad.~Publ., Dordrecht, 1991).

\bibitem{ReidSmith} Colin D.~Reid and Simon M.~Smith (with an appendix by Stephan Tornier),
	Groups acting on trees with Tits’ independence property (P).
	{\it Preprint}, (2022) arXiv:2002.11766v2.

\bibitem{Serre:trees}
	Jean-Pierre Serre,
	{\it Trees}.
	Springer Monogr.~in Math.~(Springer, Berlin, 2003).
	
\bibitem{SamShep}
    Sam Shepherd (with appendix by Giles Gardam and Daniel J.~Woodhouse),
    Two Generalisations of Leighton's Theorem.
    {\it Preprint}, (2019) to appear in {\it Groups, Geometry, and Dynamics.} arXiv:1908.00830.

\bibitem{SmithDigraph}
    Simon M.~Smith,
    Infinite primitive directed graphs.
    {\it J.~Algebr.~Comb.}
    31 (2010), 131--141.

\bibitem{SmithSubdegrees}
    Simon M.~Smith,
    Subdegree growth rates of infinite primitive permutation groups.
    {\it J.~Lond.~Math.~Soc.} 82 (2010), 526--548.

\bibitem{SmithDuke}
    Simon M.~Smith, A product for permutation groups and topological groups,
    {\it Duke Math.~J.} {166} (2017), 2965--2999.

\bibitem{SmithPrimitive}
    Simon M.~Smith, The structure of primitive permutation groups with finite suborbits and t.d.l.c.~groups admitting a compact open subgroup that is maximal.
    {\it Preprint}, (2019) arXiv:1910.13624v2.

\bibitem{ThomassenWoess}
    Carsten Thomassen and Wolfgang Woess,
    Vertex-transitive graphs and accessibility.
    {\it J.~Comb.~Theory, Ser.~B} {58} (1993), 248--268.

\bibitem{Tornier}
    Stephan Tornier, 
    Groups acting on trees with prescribed local action.
    {\it Preprint}, (2020) arXiv:2002.09876v3.

\bibitem{Tits70}
    Jacques Tits, Sur le groupe des automorphismes d'un arbre.
    In: {\it Essays on topology and related topics (M\'emoires d\'edi\'es
      \`a Georges de Rham)}, (Springer, New York, 1970), pp.~188--211.

\bibitem{weiss78}
    Richard Weiss,
    s-Transitive graphs,
    {\em Algebraic methods in graph theory} {\bf 25} (1978), 827--847.
\end{thebibliography}
\end{document}